\documentclass[a4paper, 11pt]{article}
\usepackage{graphicx}
\usepackage{pdfpages}
\usepackage{framed}
\usepackage[normalem]{ulem}
\usepackage{amsmath}
\usepackage{amsthm}
\usepackage{amssymb}
\usepackage{amsfonts}
\usepackage{bm}
\usepackage{mathtools}
\usepackage{enumerate}
\usepackage{enumitem}
\usepackage{pdflscape}
\usepackage{caption}
\usepackage[nottoc]{tocbibind}
\usepackage[title]{appendix}
\usepackage{afterpage,natbib,lipsum}
\usepackage[unicode, psdextra, linkcolor=blue, citecolor=blue, colorlinks=true, pdfencoding=auto]{hyperref} 
\usepackage{subcaption}
\usepackage[utf8]{inputenc}
\usepackage[top = 3.5cm, bottom = 3.5cm, left = 3cm, right = 3cm]{geometry}
\usepackage{array, makecell, siunitx}
\usepackage{lscape}
\usepackage{empheq}
\usepackage{accents}
\usepackage{xcolor}
\usepackage[ruled,vlined]{algorithm2e}
\usepackage{authblk}
\usepackage{array}
\usepackage{multirow}
\usepackage{threeparttable, tabularx}
\usepackage{bookmark, xpatch}
\usepackage{textgreek}

\numberwithin{equation}{section}

\newcommand{\dfF}{\mathrm{df}_{\rm F}}
\newcommand{\dfR}{\mathrm{df}_{\rm R}}
\newcommand{\bx}{\mathbf{x}}
\newcommand{\bX}{\mathbf{X}}
\newcommand{\by}{\mathbf{y}}
\newcommand{\bI}{\mathbf{I}}
\newcommand{\bH}{\mathbf{H}}
\newcommand{\bh}{\mathbf{h}}
\newcommand{\bW}{\mathbf{W}}

\newcommand{\bA}{\mathbf{A}}

\newcommand{\bs}{\mathbf{s}}
\newcommand{\bT}{\mathbf{T}}
\newcommand{\bz}{\mathbf{z}}
\newcommand{\bZ}{\mathbf{Z}}
\newcommand{\bU}{\mathbf{U}}
\newcommand{\bu}{\mathbf{u}}
\newcommand{\bV}{\mathbf{V}}

\newcommand{\bF}{\mathbf{F}}
\newcommand{\bG}{\mathbf{G}}
\newcommand{\bg}{\mathbf{g}}
\newcommand{\bR}{\mathbf{R}}
\newcommand{\br}{\mathbf{r}}

\newcommand{\bQ}{\mathbf{Q}}
\newcommand{\bq}{\mathbf{q}}

\newcommand{\bM}{\mathbf{M}}

\newcommand{\cS}{\mathcal{S}}

\newcommand{\bmu}{\bm{\mu}}
\newcommand{\bbeta}{\bm{\beta}}

\newcommand{\bSigma}{\bm{\Sigma}}

\newcommand{\bphi}{\bm{\phi}}
\newcommand{\ErrF}{\mathrm{ErrF}}
\newcommand{\ErrR}{\mathrm{ErrR}}
\newcommand{\ErrT}{\mathrm{ErrT}}
\newcommand{\OptF}{\mathrm{OptF}}
\newcommand{\OptR}{\mathrm{OptR}}
\newcommand{\wErrF}{\mathrm{wErrF}}
\newcommand{\wErrR}{\mathrm{wErrR}}
\newcommand{\wErrT}{\mathrm{wErrT}}
\newcommand{\wOptF}{\mathrm{wOptF}}
\newcommand{\wOptR}{\mathrm{wOptR}}

\newcommand{\trace}{{\rm tr}}

\newcommand{\R}{\mathbb{R}}

\newcommand{\E}{\mathrm{E}}
\newcommand{\T}{\top}
\newcommand{\var}{\mathrm{Var}}
\newcommand{\cov}{\mathrm{Cov}}
\newcommand{\prob}{\mathrm{P}}

\newtheorem{theorem}{Theorem}[section]

\newtheorem{proposition}{Proposition}[section]

\theoremstyle{definition}

\theoremstyle{definition}

\theoremstyle{definition}
\newtheorem{remark}{Remark}[section]

\title{On Measuring Model Complexity in Heteroscedastic\\ Linear Regression}
\author{Bo Luan}
\author{Yoonkyung Lee}
\author{Yunzhang Zhu}
\affil{The Ohio State University}
\date{}

\begin{document}
\maketitle

\begin{abstract}
	Heteroscedasticity is common in real world applications and is often handled by incorporating case weights into a modeling procedure. Intuitively, models fitted with different weight schemes would have a different level of complexity depending on how well the weights match the inverse of error variances. However, existing statistical theories on model complexity, also known as model degrees of freedom, were primarily established under the assumption of equal error variances. In this work, we focus on linear regression procedures and seek to extend the existing measures to a heteroscedastic setting. Our analysis of the weighted least squares method reveals some interesting properties of the extended measures. In particular, we find that they depend on both the weights used for model fitting and those for model evaluation. Moreover, modeling heteroscedastic data with optimal weights generally results in fewer degrees of freedom than with equal weights, and the size of reduction depends on the unevenness of error variance. This provides additional insights into weighted modeling procedures that are useful in risk estimation and model selection.
	\vspace{8pt}
	
\end{abstract}

\section{Introduction}

Heteroscedasticity is common in real world data and needs careful consideration for modeling. Mishandling it may lead to inefficient estimates for model parameters and unreliable predictions. To deal with heteroscedasticity, a commonly used strategy is to incorporate case weights into a modeling procedure. If the error variances associated with cases were known, ideally we would specify the weights inversely proportional to the error variances according to the normal likelihood. However, such case-specific variability is usually unknown a priori, and weights need to be estimated. This generally results in various weight schemes depending on how they are estimated. For models with different weight schemes, we are interested in whether and how their complexity measures vary. A deeper understanding of this issue will offer insights into the impact of weight schemes on model performance in relation to the true error variance and further provide guidance on weight selection and ensuing model selection. 

In statistics, the term ``model degrees of freedom'' is used as a measure of the complexity or flexibility of a model \citep{tibshirani1987local, ye1998measuring, efron2004estimation}. For least squares linear regression as an example, the model degrees of freedom is defined as the number of free parameters in the model, which coincides with the number of predictors. More generally, the degrees of freedom of a model can be defined through the expected downward bias of the training error as an estimate of the prediction error. This expected bias is also known as the \textit{expected optimism} \citep{efron1986biased}. Intuitively, models fitted with different weight schemes would differ in complexity, since both training and prediction errors are affected by case weights. However, to the best of our knowledge, no complexity measure has been formally established for weighted modeling procedures. 

In this paper, we examine the expected optimism for in-sample prediction as well as out-of-sample prediction with case weights. The classical optimism theory \citep{efron1986biased} takes in-sample prediction error as the primary target. For a standard regression setting with constant error variance, the model degrees of freedom is defined through the covariances between observed responses and their fitted values. A shortcoming of this notion of model degrees of freedom is that it fails to account for the flexibility of a model at new feature values, which could be a major issue with high dimensional data. To overcome this issue, \cite{luan2021predictive} proposed a different version of model degrees of freedom, called the \textit{predictive model degrees of freedom}, based on the out-of-sample prediction error. Both versions of model degrees of freedom have been derived under the assumption of equal error variance for regression models, and they do not directly apply to modeling procedures with case weights.

We aim to extend the optimism theory and the two notions of model degrees of freedom to a more general setting with heteroscedastic data for linear regression procedures. First, we assume that the error variance is a function of covariates rather than constant. This extension has also been considered in \cite{anscombe1961examination}, \cite{bickel1978using}, \cite{box1974correcting}, \cite{carroll1981robust}, for example. Second, unlike the existing theories that evaluate a model with squared loss, we use a weighted squared loss as an evaluation metric for a given weight scheme. Using these \textit{weights for evaluation}, we define analogous expected optimism and model degrees of freedom.

As an example of linear regression procedures, we examine the two extended measures of model degrees of freedom for the weighted least squares method in two scenarios. In the first scenario of subset regression, we study their monotonicity as the subset size increases. According to these two measures, the weighted least squares model has smaller degrees of freedom than the ordinary least squares model of the same size, and the size of reduction in the degrees of freedom depends on the unevenness of the error variances. In the second scenario, we fix the set of covariates as well as weights for evaluation and examine the impact of the weight scheme for fitting a weighted least squares model to data, which we call \textit{weights for model fitting}. From sensitivity analyses of degrees of freedom around the optimal weights, we gain some useful insights into the relationships between the two model complexity measures and case weights for model fitting.

Using the two extended measures of model degrees of freedom, we propose estimators for both in-sample and out-of-sample prediction errors in a heteroscedastic setting. For the former, our estimator is a direct extension of Mallows's $C_p$ \citep{mallows1973some}, whereas for the latter, the proposed estimator is an adaptation of the risk estimator developed in \cite{luan2021predictive} to heteroscedastic data. In practice, the predictive model degrees of freedom is usually unknown due to its dependency on the underlying data distribution. We thereby propose methods that can estimate the predictive model degrees of freedom efficiently and use them for estimation of the out-of-sample prediction error.

These proposed estimators are then assessed through a series of simulation studies and real data analysis. First, we evaluate our out-of-sample risk estimator in terms of prediction error estimation and model selection. We find that it has smaller variance than the leave-one-out and five-fold cross validation strategies. In addition, it tends to favor slightly simpler models without compromising the predictive performance too much. Second, we evaluate the impact of weight misspecification on model selection for both in-sample and out-of-sample predictions. In general, ignoring heteroscedasticity may lead to misidenfication of the model size and inflation in prediction risk. 

This paper is organized as follows. In Section \ref{sec: review_df}, we first review the two notions of model degrees of freedom defined under the assumption of equal error variance. We then propose their extended versions in a heteroscedastic setting in Section \ref{sec: extend_df}. We examine properties of the extended measures for the weighted least squares method in Section \ref{sec: wls}, and propose approaches to estimation of the predictive model degrees of freedom and out-of-sample prediction error in Section \ref{sec: pred_err_est}. In Section \ref{sec: numerical_study}, we evaluate the proposed estimators through simulation studies and real data analysis. We then conclude the paper in Section \ref{sec: conclusion}. For conciseness, all proofs will be deferred to Appendix.

\section{Review of Model Complexity Theories} \label{sec: review_df}

In this section, we review the existing theories for the complexity of regression models. To start with, we first introduce the notation that will be used in this work. Consider a regression setting with a continuous response $Y \in \R$. For a given set of covariates $X = (X_1,\ldots,X_d)$, assume that $Y = \mu(X;\bbeta) + \varepsilon$, where $\mu$ is the mean regression function with unknown parameter $\bbeta$, and $\varepsilon$ is some random error such that $\E(\varepsilon \vert X) = 0$ and $\var(\varepsilon \vert X) = \sigma_\varepsilon^2 > 0$. Here, we allow $X_1,\ldots,X_d$ to be either numerical or categorical. 

Assume that a dataset of size $n$ is used for model fitting. Let $\by = (y_1,\ldots,y_n)^\T$ be the vector of responses and $\bX = (\bx_1,\ldots,\bx_n)^\T \in \R^{n \times m}$ be the full design matrix that includes all $d$ covariates. We encode a categorical variable with $L$ levels by $L-1$ binary columns in $\bX$. Thus, $m > d$ if there is at least one categorical variable with more than two levels. 

For brevity, we abbreviate $\mu(\bx_i;\bbeta)$ by $\mu_i$ and write $\bmu = (\mu_1,\ldots,\mu_n)^\T$. For a given modeling procedure, define $\hat{\mu}_i$ to be the fitted value of $y_i$ and write $\hat{\bmu} = (\hat{\mu}_i, \ldots, \hat{\mu}_n)^\T$. With the training data $(\bX, \by)$, $\hat{\bmu}$ is obtained by minimizing the (unweighted) training error
$$
	\ErrT_{\bX,\by} = \frac{1}{n}\sum_{i=1}^n (y_i - \hat{\mu}_i)^2 = \frac{1}{n}\Vert \by - \hat{\bmu} \Vert^2.
$$

\subsection{Classical Model Degrees of Freedom} \label{subsec: review_classical_df}

We first review the classical optimism theory presented in \cite{efron1986biased,efron2004estimation} that focuses on in-sample prediction error. Let $\tilde{y}_i$ be an independent copy of $y_i$ given $\bx_i$. The (unweighted) in-sample prediction error is defined as
$$
	\ErrF_{\bX,\by} = \frac{1}{n}\E(\Vert \tilde{\by} - \hat{\bmu} \Vert^2 \vert \bX, \by),
$$
where the expectation is with respect to $\tilde{\by} = (\tilde{y}_1,\ldots,\tilde{y}_n)^\T$. 

The expected optimism, denoted by $\OptF_\bX$, is the expected difference between the prediction error $\ErrF_{\bX,\by}$ and training error $\ErrT_{\bX,\by}$. The expectation is taken with respect to $\by$ given $\bX$. \cite{efron1986biased} showed that
\begin{equation}\label{eq: OptF_unweighted}
	\OptF_{\bX} \coloneqq \E(\ErrF_{\bX,\by} - \ErrT_{\bX,\by}\vert \bX) = \frac{2}{n} \sum_{i=1}^n \cov(y_i, \hat{\mu}_i\vert \bX).
\end{equation}
Following this covariance representation of the expected optimism, \cite{efron2004estimation} then formally defined the model degrees of freedom of a fitted model $\hat{\bmu}$ as
\begin{equation} \label{eq: dfF_unweighted}
    \dfF = \sum_{i=1}^n \frac{\cov(y_i, \hat{\mu}_i\vert \bX)}{\sigma_\varepsilon^2}.
\end{equation}
In particular, for a linear modeling procedure $\hat{\bmu} = \bH \by$, where $\bH = (h_{ij}) \in \R^{n \times n}$ depends only on $\bX$, it is easy to show that $\cov(y_i, \hat{\mu}_i\vert \bX) = \sigma_\varepsilon^2 h_{ii}$ and $\dfF = \trace(\bH)$. Here, $\bH$ is usually called the \textit{hat matrix}.

When the input space is large, the in-sample prediction error may be limited in measuring how well a modeling procedure will generalize at new feature values that arise in the future as in many prediction-oriented applications. As a result, the classical model degrees of freedom $\dfF$ needs some adjustment when prediction beyond the training data is concerned. This motivates the use of out-of-sample prediction error for an alternative definition of model degrees of freedom, which we will discuss next.

\subsection{Predictive Model Degrees of Freedom} \label{subsec: review_predictive_df}

In this section, we review the model complexity theory for out-of-sample prediction proposed by \cite{luan2021predictive}. Let $(\bx_\ast,\varepsilon_\ast)$ be an independent copy of $(\bx_i,\varepsilon_i)$ and define $\mu_\ast = \mu(\bx_\ast; \bbeta)$ and $y_\ast = \mu_\ast + \varepsilon_\ast$. Let $\hat{\mu}_\ast$ be a prediction of $y_\ast$. The (unweighted) out-of-sample prediction error is defined as
$$
	\ErrR_{\bX,\by} = \E[(y_\ast - \hat{\mu}_\ast)^2 \vert \bX, \by],
$$
where the expectation is taken with respect to $(\bx_\ast, y_\ast)$. Following the framework for in-sample prediction, we define the expected optimism as $\OptR_\bX = \ErrR_\bX - \ErrT_\bX$, where $\ErrR_\bX$ and $\ErrT_\bX$ are the expectations of $\ErrR_{\bX,\by}$ and $\ErrT_{\bX,\by}$ with respect to $\by$ respectively.

Consider a linear procedure with hat matrix $\bH$. For any new covariate vector $\bx_\ast$, the procedure predicts the response at $\bx_\ast$ using a linear combination of $y_i$'s, denoted by $\hat{\mu}_\ast = \bh_\ast^\T \by$, where $\bh_\ast \in \R^n$ depends only on $\bx_\ast$ and $\bX$. This coefficient vector $\bh_\ast$ is called the \textit{hat vector} at $\bx_\ast$ reminiscent of the hat matrix. Then, we can write the bias-variance decomposition for $\ErrR_\bX$ and $\ErrT_\bX$ as
\begin{align}
	& \ErrR_{\bX} = \sigma_\varepsilon^2 + \E[(\mu_\ast - \bh_\ast^\T \bmu)^2 \vert \bX] + \sigma_\varepsilon^2 \E(\Vert \bh_\ast \Vert^2 \vert \bX), \label{eq: bias_var_ErrR_unweighted}\\
	& \ErrT_\bX = \sigma_\varepsilon^2 + \frac{1}{n}\Vert \bmu - \bH\bmu \Vert^2 + \frac{1}{n} \sigma_\varepsilon^2 \trace(\bH^\T \bH - 2 \bH), \label{eq: bias_var_ErrT_unweighted}
\end{align}
where the three terms on the right hand side are the irreducible error, (squared) bias and variance respectively. The expected optimism $\OptR_\bX$ can then be expressed as
\begin{equation} \label{eq: OptR_unweighted_decomposition}
\begin{aligned}
	\OptR_\bX & = \Delta B_\bX + \frac{2}{n} \sigma_\varepsilon^2 \left[\dfF + \frac{n}{2}\left(\E(\Vert \bh_\ast\Vert^2 \vert \bX) - \frac{1}{n}\trace( \bH \bH^\T)\right)\right],
\end{aligned}
\end{equation}
where $\Delta B_\bX = \E[(\mu_\ast - \bh_\ast^\T \bmu)^2 \vert \bX] - \frac{1}{n}\Vert \bmu - \bH\bmu \Vert^2$ is the \textit{excess bias} and $\dfF$ is the classical model degrees of freedom defined in \eqref{eq: dfF_unweighted}. The second term in the right hand side of \eqref{eq: OptR_unweighted_decomposition} is the \textit{excess variance}, and it doesn't depend on the underlying true model $\mu$. Further, the quantity in the square brackets depends only on the modeling procedure through the hat matrix $\bH$ and hat vector $\bh_\ast$. Since model complexity is an intrinsic property of a modeling procedure, \cite{luan2021predictive} defined this quantity as a model complexity measure under the out-of-sample prediction setting. Denoting it by $\dfR$:
\begin{equation} \label{eq: dfR_unweighted}
	\dfR = \dfF + \frac{n}{2}\left(\E(\Vert \bh_\ast\Vert^2 \vert \bX) - \frac{1}{n}\trace(\bH^\T \bH)\right),
\end{equation}
they call $\dfR$ the \textit{predictive model degrees of freedom} to differentiate it from the classical model degrees of freedom and emphasize its use in an out-of-sample prediction setting.

\section{Model Complexity in a Heteroscedastic Setting} \label{sec: extend_df}

In this section, we extend the existing theories on model complexity to a heteroscedastic setting. Our extensions include two aspects. First, we allow heteroscedasticity in error $\varepsilon$. In particular, we assume that $\var(\varepsilon \vert X) = \sigma_\varepsilon^2 \tau(X)$, where $\tau$ is a function satisfying $\tau(\cdot) > 0$ and $\E[\tau(X)] = 1$. We write $\tau_i = \tau(\bx_i)$ for simplicity. Second, for model evaluation, we consider a general weighted squared error loss. Let $w \colon \R^m \to \R^+$ be a deterministic function and write $w_i = w(\bx_i)$. We define the weighted training error given $(\bX,\by)$ and $w$ as
$$
	\wErrT_{\bX,\by} = \frac{1}{n}\sum_{i=1}^n w_i (y_i - \hat{\mu}_i)^2 = \frac{1}{n}\Vert \by - \hat{\bmu} \Vert_\bW^2,
$$
where $\bW = \mathrm{diag}(w_1,\ldots,w_n)$ and $\Vert \mathbf{a} \Vert_\bW^2 = \mathbf{a}^\T \bW \mathbf{a}$ for $\mathbf{a} \in \R^n$. In particular, if we assume $y_i \sim \mathcal{N}(\mu_i, \sigma_\varepsilon^2 \tau_i)$ given $\bx_i$ and take $w_i = 1/\tau_i$, then minimizing the weighted training error is equivalent to maximizing the log-likelihood of $\by$ given $\bX$. For in-sample and out-of-sample prediction errors, we incorporate case weights in a similar way and define
$$
    \wErrF_{\bX,\by} = \frac{1}{n}\E(\Vert \tilde{\by} - \hat{\bmu} \Vert_\bW^2 \vert \bX, \by) \text{ and } \wErrR_{\bX,\by} = \E[w_\ast (y_\ast - \hat{\mu}_\ast)^2 \vert \bX, \by],
$$
where $\tilde{\by}$ is an independent copy of $\by$ given $\bX$ and $w_\ast = w(\bx_\ast)$.
 
\subsection{Classical Model Degrees of Freedom}
We first look at the classical model degrees of freedom. With the weighted training error and in-sample prediction error, it is easy to adjust the expected optimism \eqref{eq: OptF_unweighted} with weights as
$$
	\wOptF_{\bX} = \frac{2}{n} \sum_{i=1}^n w_i \cov(y_i, \hat{\mu}_i\vert \bX).
$$
Based on the expected optimism, we define the extended model degrees of freedom as
\begin{equation} \label{eq: dfF_weighted}
    \dfF = \sum_{i=1}^n \frac{w_i}{\bar{w}} \frac{\cov(y_i, \hat{\mu}_i\vert \bX)}{\sigma_\varepsilon^2},
\end{equation}
where $\bar{w} = \sum_{i=1}^n w_i / n$. The use of $w_i/\bar{w}$ as a normalized weight for each covariance summand guarantees that $\dfF$ is invariant to scaling of the weight. In addition, the normalizing constant $\bar{w}$ ensures that the normalized weights sum up to the sample size $n$. In particular, when $w_1 = \cdots = w_n$, the extended definition just reduces to the model degrees of freedom \eqref{eq: dfF_unweighted} for equal error variances.

Using the extended model degrees of freedom $\dfF$, we can immediately get an unbiased estimator of $\wErrF_{\bX}$, which is given by
$$
    \widehat{\wErrF} = \wErrT_{\bX, \by} + \frac{2}{n} \bar{w} \sigma_\varepsilon^2 \dfF.
$$

For a linear procedure with hat matrix $\bH$, we have $\cov(y_i, \hat{\mu}_i\vert \bX) = h_{ii} \sigma_\varepsilon^2 \tau_i$. Then \eqref{eq: dfF_weighted} is simplified to
\begin{equation} \label{eq: dfF_linear_weighted}
    \dfF = \sum_{i=1}^n \frac{w_i}{\bar{w}} h_{ii} \tau_i = \frac{1}{\bar{w}} \trace(\bH \bT \bW),
\end{equation}
where $\bT = \mathrm{diag}(\tau_1,\ldots,\tau_n)$. When $\bW = \bT = \bI_n$, this further reduces to $\dfF = \trace(\bH)$ in the homoscedastic setting. The following remarks demonstrate some interesting properties of this extension \eqref{eq: dfF_linear_weighted}.

\begin{remark}[Optimal weights] \label{rmk: optimal_weights}
    Assume the optimal weights are used for model evaluation, i.e., $w(\cdot) \propto 1/\tau(\cdot)$. Then we have
    \begin{equation} \label{eq: dfF_weighted_optimal}
        \dfF = \zeta_\tau \trace(\bH),
    \end{equation}
    where $\zeta_\tau = \left(\frac{1}{n}\sum_{i=1}^n \tau_i^{-1}\right)^{-1}$ is the harmonic mean of $\tau_1,\ldots,\tau_n$. Since $\tau_i$'s are all positive, their harmonic mean is always no greater than the arithmetic mean $\bar{\tau}$, and the equality holds if and only if all $\tau_i$'s are equal. This implies that
    $$
        \dfF = \zeta_\tau \trace(\bH) \leq \bar{\tau} \trace(\bH) \approx \trace(\bH),
    $$
    where the approximation in the last step is due to $\bar{\tau} \approx \E[\tau(\bx)] = 1$. Note that $\trace(\bH)$ is the classical model degrees of freedom of a linear procedure under the homoscedastic setting. This suggests that the unevenness in error variances leads to a deflation in model complexity. 

\end{remark}

\begin{remark}[Equal error variances]
    When error variances are all equal, i.e., $\tau_1 = \cdots = \tau_n = 1$, \eqref{eq: dfF_linear_weighted} becomes $\dfF = \sum_{i=1}^n (w_i/\bar{w}) h_{ii}$, which is a weighted average of the leverages $h_{ii}$. Assume $h_{ii} > 0$ for $i = 1,\ldots,n$. Then, for any $w_1,\ldots,w_n > 0$,  we have $n \min_{i} h_{ii} \leq \dfF \leq n \max_{i} h_{ii}$.
    In other words, $\dfF$ is maximized (minimized) when all weights are assigned to the case with the largest (smallest) leverage. Since leverages are directly related to the impact of observations on the model, putting more weights on more (less) influential observations will make the model more (less) flexible. The extended model degrees of freedom properly reflects this relationship.
\end{remark}

\begin{remark}[Duplicate data]
    Assume $\bx_1,\ldots,\bx_k$ are $k$ distinct covariate vectors in the training data with $n_1,\ldots,n_k$ replicates respectively ($n = \sum_{i=1}^k n_i$). For each $i = 1,\ldots,k$, let $y_{i,1},\ldots y_{i,n_i}$ be i.i.d.\ observations given $\bx_i$ with constant variance $\sigma_\varepsilon^2$, and define $\bar{y}_i = \sum_{j=1}^{n_i} y_{i,j} / n_i$. Note that
    $$
        \sum_{i=1}^k \sum_{j=1}^{n_i} (y_{i,j} - \hat{\mu}_i)^2 = \sum_{i=1}^k \sum_{j=1}^{n_i} (y_{i,j} - \bar{y}_i)^2 + \sum_{i=1}^k n_i (\bar{y}_i - \hat{\mu}_i)^2,
    $$
    and the first term on the right hand side of the above equation doesn't depend on the model. Thus, a weighted least squares model based on the collapsed data $\lbrace (\bx_i, \bar{y}_i) \vert i = 1,\ldots,k \rbrace$ with weights $w_i = n_i$ is identical to the ordinary least squares model based on the full data $\lbrace (\bx_i, y_{i,j}) \vert i = 1,\ldots,k;\, j = 1,\ldots, n_i \rbrace$. This dual representation of the model brings a question about its degrees of freedom. If the collapsed data were to be viewed as a collection of genuine $k$ independent observations with unequal weights, Remark \ref{rmk: optimal_weights} appears to suggest that the weighted model with a linear procedure has fewer degrees of freedom than the unweighted model fit to the same $k$ observations. However, this difference in model complexity comes from the scale of weights tied to the sample size we assume. If the weights are scaled to match the original sample size $n$ rather than the distinct number of observations $k$, the dual representation of the model fit to duplicate data leads to the same model degrees of freedom.
\end{remark}

\subsection{Predictive Model Degrees of Freedom}
We now extend the notion of model complexity for the heteroscedastic setting to the out-of-sample prediction framework for all linear procedures. Let $w_\ast = w(\bx_\ast)$ and $\tau_\ast = \tau(\bx_\ast)$, and define the weighted out-of-sample prediction error as
$$
	\wErrR_{\bX,\by} = \E[w_\ast (y_\ast - \hat{\mu}_\ast)^2 \vert \bX, \by],
$$
where the expectation is with respect to $(\bx_\ast, y_\ast)$. Consider a linear procedure with hat matrix $\bH$ and hat vector $\bh_\ast$ at $\bx_\ast$. Similar to \eqref{eq: bias_var_ErrR_unweighted} and \eqref{eq: bias_var_ErrT_unweighted} in the homoscedastic setting, we have the following bias-variance decompositions for $\wErrR_\bX$ and $\wErrT_\bX$: 
\begin{align}
	& \wErrR_{\bX} = \sigma_\varepsilon^2 \E(w_\ast \tau_\ast) + \E[w_\ast (\mu_\ast - \bh_\ast^\T \bmu)^2 \vert \bX] + \sigma_\varepsilon^2 \E(w_\ast \Vert \bh_\ast \Vert_\bT^2 \vert \bX),\label{eq: bias_var_ErrR_weighted}\\
	& \wErrT_\bX = \frac{1}{n}\sigma_\varepsilon^2\sum_{i=1}^n w_i \tau_i + \frac{1}{n}\Vert \bmu - \bH\bmu \Vert_\bW^2 + \frac{1}{n} \sigma_\varepsilon^2 \trace(\bH^\T \bW \bH \bT - 2 \bH \bT \bW), \label{eq: bias_var_ErrT_weighted}
\end{align}
where $\wErrR_{\bX}$ is the expectation of $\wErrR_{\bX,\by}$ with respect to $\by$. See Appendix \ref{pf: bias_var_Err_weighted} for detailed derivation of \eqref{eq: bias_var_ErrR_weighted} and \eqref{eq: bias_var_ErrT_weighted}. As a consequence, the expected optimism $\wOptR_\bX \coloneqq \wErrR_\bX - \wErrT_\bX$ is of the form
\begin{equation} \label{eq: OptR_weighted_decomposition}
\begin{aligned}
	\wOptR_\bX 
	& = \sigma_\varepsilon^2\left(\E(w_\ast \tau_\ast) - \frac{1}{n}\sum_{i=1}^n w_i \tau_i \right)\\
	& \qquad + \left(\E[w_\ast (\mu_\ast - \bh_\ast^\T \bmu)^2 \vert \bX] - \frac{1}{n}\Vert \bmu - \bH\bmu \Vert_\bW^2\right)\\ 
	& \qquad + \frac{2}{n} \bar{w} \sigma_\varepsilon^2 \left[\dfF + \frac{n}{2\bar{w}}\left(\E(w_\ast \Vert \bh_\ast\Vert_\bT^2 \vert \bX) - \frac{1}{n}\trace(\bW \bH \bT \bH^\T)\right)\right],
\end{aligned}
\end{equation}
where $\dfF =\trace(\bH \bT \bW)/\bar{w}$ is the model degrees of freedom for a linear procedure defined in \eqref{eq: dfF_linear_weighted}. We call the first term on the right the \textit{excess irreducible error}, which equals $0$ in the homoscedastic setting due to $\tau$ and $w$ both being constant. In general, this term should be negligible as long as the sample size is large. The second and third terms are respectively the excess bias and excess variance adjusted by the weight function $w$. Following \cite{luan2021predictive}, we still use $\Delta B_\bX$ to denote the excess bias. We then propose the extended predictive model degrees of freedom based on the excess variance. Formally, we define it as
\begin{equation} \label{eq: dfR_weighted}
	\dfR = \dfF + \frac{n}{2\bar{w}}\left(\E(w_\ast \Vert \bh_\ast\Vert_\bT^2 \vert \bX) - \frac{1}{n}\trace(\bW \bH \bT \bH^\T)\right).
\end{equation}
In particular, when $w(\cdot)$ is constant, we have
$$
	\dfR = \dfF + \frac{n}{2}\left(\E(\Vert \bh_\ast\Vert_\bT^2 \vert \bX) - \frac{1}{n}\trace(\bH \bT \bH^\T)\right).
$$
If further $\tau_i$'s are all equal, this reduces to \eqref{eq: dfR_unweighted}.

The expression in \eqref{eq: dfR_weighted} suggests that the predictive model degrees of freedom $\dfR$ adjusts the classical model degrees of freedom $\dfF$ with an extra term that accounts for out-of-sample prediction. In fact, this adjustment can be expressed and interpreted in two different ways, which we summarize in the following proposition.
\begin{proposition}\label{prop: dfR_cov_representation}
	For a linear procedure, the adjustment $\dfR - \dfF$ can be expressed as
	\begin{enumerate}
	    \item[(i)] Covariance penalty representation
	    $$
		    \dfR - \dfF = \frac{n}{2\bar{w}}\sum_{i=1}^n \tau_i \left[\E\left(w_\ast \frac{\cov^2(y_i, \hat{\mu}_\ast \vert \bx_\ast, \bX)}{(\sigma_\varepsilon^2 \tau_i)^2} \Bigg\vert \bX\right) - \frac{1}{n} \sum_{j=1}^n w_j \frac{\cov^2(y_i, \hat{\mu}_j \vert \bX)}{(\sigma_\varepsilon^2 \tau_i)^2}\right];
	    $$
	    
	    \item[(ii)] Generalized degrees of freedom (GDF) representation
	    $$
	    	\dfR - \dfF = \frac{n}{2\bar{w}}\sum_{i=1}^n \tau_i \left[\E\left(w_\ast \left(\frac{\partial \E(\hat{\mu}_\ast \vert \bx_\ast, \bX) }{\partial \mu_i}\right)^2 \Bigg\vert \bX\right) - \sum_{j=1}^n w_j \left(\frac{\partial \E(\hat{\mu}_j \vert \bX) }{\partial \mu_i}\right)^2\right].
	    $$
	\end{enumerate}
\end{proposition}

The covariance penalty representation unifies the two measures of model degrees of freedom $\dfF$ and $\dfR$ by quantifying the level of dependency of the predicted values on the observed values. More specifically, the adjustment $\dfR - \dfF$ compares the overall level of such dependency for out-of-sample predictions against that for in-sample ones. The GDF representation gets its name from the generalized degrees of freedom (GDF) proposed by \cite{ye1998measuring}. In his work, he defined $\sum_{i=1}^n \frac{\partial \E(\hat{\mu}_i \vert \bX) }{\partial \mu_i}$ as the GDF for a modeling procedure $\hat{\bmu}$. The GDF representation of $\dfR - \dfF$ expresses this adjustment in terms of the rate of change of the predicted values on the test data against that on the training data with respect to the true mean components $\mu_i$'s.

\section{Weighted Least Squares} \label{sec: wls}

We study the two proposed measures of model degrees of freedom for the weighted least squares method in this section. Assume $\bX \in \R^{n \times m}$ has full column rank. Let $q \colon \R^m \to \R_+$ be a known weight function, and define $q_i = q(\bx_i)$ and $\bQ = \mathrm{diag}(q_1, \ldots, q_n)$. For the weighted least squares problem of $\min_{\mathbf{b} \in \R^p} \Vert \by - \bX \mathbf{b} \Vert_\bQ^2$, the solution has the explicit form of $\hat{\bbeta} = (\bX^\T \bQ \bX)^{-1} \bX^\T \bQ \by$. The corresponding hat matrix and hat vector are then given by $\bH = \bX (\bX^\T \bQ \bX)^{-1} \bX^\T \bQ $ and $\bh_\ast = \bQ \bX (\bX^\T \bQ \bX)^{-1} \bx_\ast$.

Consider evaluating the model using a weighted squared error loss with weight function $w$. Let $\bSigma = \E(w_\ast \bx_\ast \bx_\ast^\T) / \bar{w}$. Then it is easy to show that
\begin{align}
    \dfF & = \frac{1}{\bar{w}} \trace[(\bX^\T \bQ \bX)^{-1} \bX^\T \bQ \bT \bW \bX], \label{eq: dfF_wls_general}\\
    \dfR & = \dfF + \frac{n}{2} \trace\left[(\bX^\T \bQ \bX)^{-1} \bX^\T \bQ \bT \bQ \bX (\bX^\T \bQ \bX)^{-1} \left(\bSigma - \frac{1}{n\bar{w}}\bX^\T \bW \bX\right)\right]. \label{eq: dfR_wls_general}
\end{align}

\subsection{Subset Regression} \label{subsec: subset_regression} 
We first study $\dfF$ and $\dfR$ of the weighted least squares model in the subset regression setting. For simplicity, we assume the weight scheme used for modeling fitting is the same as that for model evaluation, i.e., $q = w$. Let $\cS_1 \subset \cdots \subset \cS_d \coloneqq \lbrace 1,\ldots,d \rbrace$ be a sequence of nested subsets of variables such that $\vert \cS_p \vert = p$ for each $p = 1,\ldots,d$. Let $\bX_p$ and $\bSigma_p$ be the submatrices of $\bX$ and $\bSigma$ corresponding to $\cS_p$ respectively. Consider regressing $\by$ on $\bX_p$ using the weighted least squares method with weight function $w$. Then \eqref{eq: dfF_wls_general} and \eqref{eq: dfR_wls_general} are simplified to
\begin{align*}
    \dfF(p) & = \frac{1}{\bar{w}} \trace[(\bX_p^\T \bW \bX_p)^{-1} \bX_p^\T \bW \bT \bW \bX_p], \\
    \dfR(p) & = \frac{1}{2} \dfF(p) + \frac{n}{2} \trace[(\bX_p^\T \bW \bX_p)^{-1} \bX_p^\T \bW \bT \bW \bX_p (\bX_p^\T \bW \bX_p)^{-1} \bSigma_p].
\end{align*}

For each $j = 1,\ldots,d$, let $\bx_{(j)}$ be the column(s) in $\bX$ associated with variable $j$. If the $j$th variable is continuous or binary, $\bx_{(j)} \in \R^n$. If the variable is categorical with $k \geq 3$ levels, $\bx_{(j)} \in \R^{n \times (k-1)}$. For the sake of convenience and conciseness, we assume $\bx_{(j)} \in \R^n$ in the rest of this section, though all results can be tailored to categorical variables with more than two levels accordingly. 

We examine the size of change in the model degrees of freedom as a covariate is added. Normalizing observations with weights, let $\tilde{\bx}_{(p+1)} = \bW^{\frac{1}{2}} \bx_{(p+1)}/\sqrt{\bar{w}}$ and $\tilde{\bX}_p = \bW^{\frac{1}{2}} \bX_p/\sqrt{\bar{w}}$. Similarly, let $\tilde{\by} = \bW^{\frac{1}{2}} \by /\sqrt{\bar{w}}$ and note that $\var(\tilde{\by} \vert \bX) = \sigma_\varepsilon^2 \bW^{\frac{1}{2}} \bT \bW^{\frac{1}{2}}$. Define $\tilde{\bH} = \tilde{\bX}_p (\tilde{\bX}_p^\T  \tilde{\bX}_p)^{-1} \tilde{\bX}_p^\T$ and $\tilde{\br}_{p+1} = (\bI_n - \tilde{\bH}) \tilde{\bx}_{(p+1)}$. Here, $\tilde{\br}_{(p+1)}$ is the residual vector when the normalized new variable vector $\tilde{\bx}_{(p+1)}$ is regressed on the normalized design matrix $\tilde{\bX}_p$. The following theorem describes the increment in $\dfF$ in terms of $\tilde{\br}_{(p+1)}$ as $p$ increases by 1. 
\begin{theorem}\label{thm: DdfF_wls_general}
    Assume $p < d < n$. Then for the weighted least squares model with weight matrix $\bW$ and design matrix $\bX_p$, $\dfF$ increases in $p$, and
    $$
        \Delta\dfF(p) \coloneqq \dfF(p+1) - \dfF(p) = \frac{\Vert \tilde{\br}_{(p+1)}\Vert_{\bU}^2}{\Vert \tilde{\br}_{(p+1)}\Vert^2},
    $$
    where $\bU = \bW^{\frac{1}{2}} \bT \bW^{\frac{1}{2}} / \bar{w}$.
\end{theorem}

Theorem \ref{thm: DdfF_wls_general} indicates that the increment in $\dfF$ generally depends on the amount of collinearity between the normalized new variable and existing variables. A special case occurs when $w \propto 1/\tau$. In this scenario, we have $\dfF = \zeta_\tau p$ and $\Delta \dfF(p) = \zeta_\tau$ (the harmonic mean of $\tau_1,\ldots,\tau_n$). Hence, $\dfF$ is a linear function of $p$. Its increment $\zeta_\tau$ is constant for all $p$ and depends only on the unevenness of the true error variances.

For the predictive model degrees of freedom $\dfR$, assume that $\bSigma$ is positive definite. For $p = 1,\ldots,d-1$, write
$$
    \bSigma_{p+1} = 
    \begin{pmatrix}
        \bSigma_p & \bphi_{p,p+1}\\
        \bphi_{p,p+1}^\T & \sigma_{p+1}^2
    \end{pmatrix},\quad p = 1,\ldots, d-1.
$$
Let $a_{p,p+1} = \sigma_{p+1}^2 - \bphi_{p,p+1}^\T \bSigma_p^{-1} \bphi_{p, p+1}$. Note that $a_{p,p+1} > 0$ when $\bSigma$ is positive definite. We derive the increment $\Delta \dfR(p) = \dfR(p+1) - \dfR(p)$ in the following theorem.
\begin{theorem}\label{thm: DdfR_wls_general}
    Assume $p < d < n$. Then for the weighted least squares model with weight matrix $\bW$ and design matrix $\bX_p$,
    \begin{equation}\label{eq: Ddf_wls_general}
    \begin{aligned}
        \Delta\dfR(p) = \frac{1}{2}\Delta\dfF(p) & + \frac{n}{2} \Delta\dfF(p) \frac{(\Vert \tilde{\bx}_{(p+1)} - \tilde{\bX}_p \bSigma_p^{-1} \bphi_{p,p+1} \Vert_\mathbf{M}^2 + a_{p,p+1})}{\Vert \tilde{\br}_{(p+1)} \Vert^2} \\
        & - n \frac{\langle \tilde{\br}_{(p+1)}, \tilde{\bx}_{(p+1)} - \tilde{\bX}_p \bSigma_p^{-1} \bphi_{p,p+1} \rangle_{\bU\mathbf{M}}}{\Vert \tilde{\br}_{(p+1)} \Vert^2}
    \end{aligned}
    \end{equation}
    where $\mathbf{M} = \tilde{\bX}_p (\tilde{\bX}_p^\T  \tilde{\bX}_p)^{-1} \bSigma_p (\tilde{\bX}_p^\T \tilde{\bX}_p)^{-1} \tilde{\bX}_p^\T$, and $\langle \cdot, \cdot\rangle_\mathbf{A}$ for a matrix $\bA$ denotes the inner product with respect to $\mathbf{A}$.
\end{theorem}

The component $\tilde{\bx}_{(p+1)} - \tilde{\bX}_p \bSigma_p^{-1} \bphi_{p,p+1}$ can be seen as the residual vector $\tilde{\br}_{(p+1)}$ at the population level, in which we simply replace $\tilde{\bX}_p^\T \tilde{\bX}_p$ and $\tilde{\bX}^\T \tilde{\bx}_{(p+1)}$ with their population versions $n \bar{w} \Sigma_p$ and $n \bar{w} \phi_{p,p+1}$ respectively. In general, $\Delta\dfR(p)$ is not necessarily positive due to the third term in the right hand side of \eqref{eq: Ddf_wls_general}. However, when $w \propto 1/\tau$, this term vanishes due to the fact that $(\bI_n - \tilde{\bH})\bU\mathbf{M} = \mathbf{O}_{n \times n}$. In this case, we have
\begin{equation}\label{eq: dfR_wls_subset_special}
    \dfR(p) = \frac{\zeta_\tau}{2} (p + n \bar{w} \trace[(\bX_p^\T \bW \bX_p)^{-1} \bSigma_p)]),
\end{equation}
which is strictly increasing in $p$. The following remark gives an approximation of \eqref{eq: dfR_wls_subset_special} under the assumption that $\bx_i$'s are multivariate normal. 

\begin{remark}
    Assume $\bX \in \R^{n \times p}$. We slightly modify the setting by assuming $\tau$ is a function of a random variable $\bu$ that is independent of $\bx$. Further, assume that $\bx_\ast$ is multivariate normal with $\E(\bx_\ast) = \mathbf{0}$. Then it can be shown that
    \begin{equation}\label{eq: dfR_wls_normal_approx}
        \dfR(p) \approx \frac{p}{2} \zeta_\tau \left(1 + \frac{n_\tau}{n_\tau - p - 1}\right),
    \end{equation}
    where $n_\tau$ is the nearest integer to $\left(\sum_{i=1}^n 1/\tau_i\right)^2 / \sum_{i=1}^n (1/\tau_i)^2$. For the derivation of \eqref{eq: dfR_wls_normal_approx}, see Appendix \ref{pf: dfR_wls_normal_approx}.
\end{remark}

\begin{figure}[t]
    \centering
    \includegraphics[scale = 0.57]{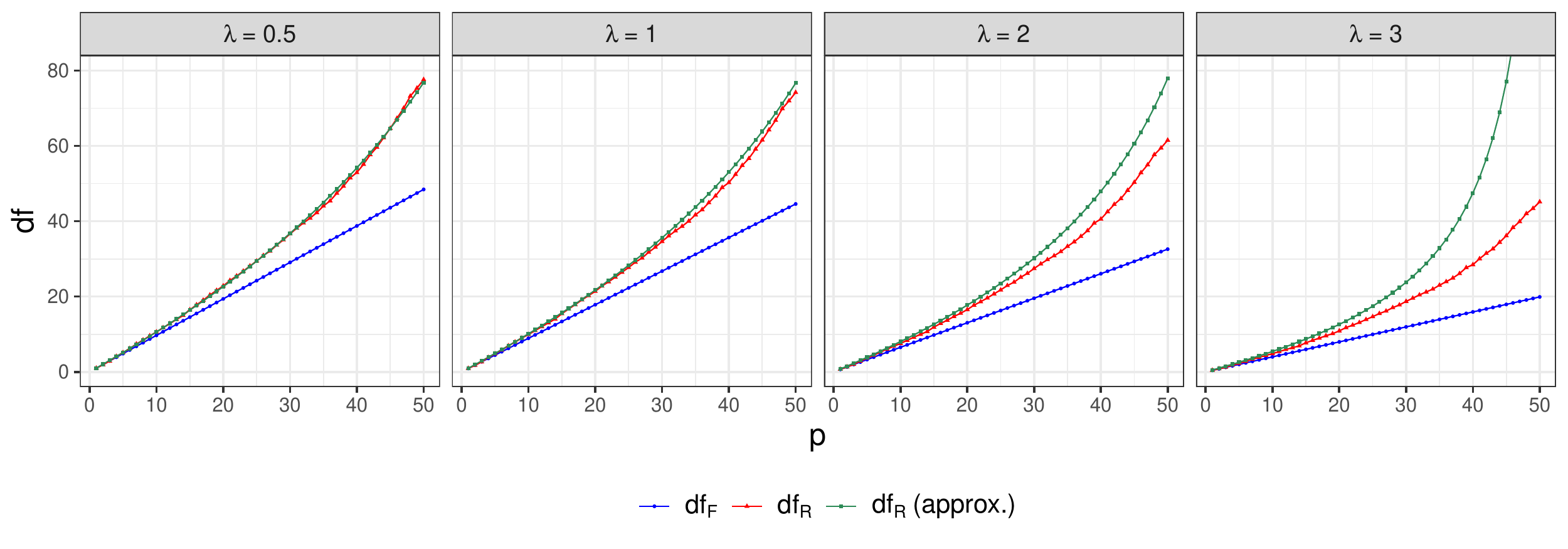}
    \caption{The approximated predictive model degrees of freedom $\dfR$ of the weighted least squares model with $w = 1/\tau$ as the weight function. The results are obtained with $\mu(\bx_i) = \sqrt{10/d} \sum_{j=1}^d x_{ij}$ and $\tau(u_i) = (1 + \vert u_i \vert)^\lambda$ ($\lambda \geq 0$), where $d=50$, $n = 100$, and $x_{ij}$'s and $u_i$ are i.i.d.\ $\mathcal{N}(0,1)$.}
    \label{fig: dfR_wls_normal_approx}
\end{figure}

Figure \ref{fig: dfR_wls_normal_approx} shows the approximation \eqref{eq: dfR_wls_normal_approx} under the assumption that $\tau(u_i) = (1 + \vert u_i \vert)^\lambda$ ($\lambda \geq 0$) and $\mu(\bx_i) = \sqrt{10/d} \sum_{j=1}^d x_{ij}$, where $x_{ij}$'s and $u_i$ are i.i.d.\ $\mathcal{N}(0,1)$. We set $n = 100$ and $d = 50$ and consider $\lambda = 0.5, 1, 2 \text{ and } 3$. The distribution of $\tau(u)$ becomes more skewed to the right as $\lambda$ increases. Consequently, the rate of increase in both $\dfF$ and $\dfR$ drops as the skewness level $\lambda$ increases. The approximation \eqref{eq: dfR_wls_normal_approx} works well when the distribution of $\tau$ is less skewed, but tends to overestimate the true predictive model degrees of freedom with more uneven $\tau_i$'s at a large $p$ due to the quickly shrinking denominator $n_\tau - p - 1$.

\begin{figure}[t]
    \centering
    \includegraphics[scale = 0.53]{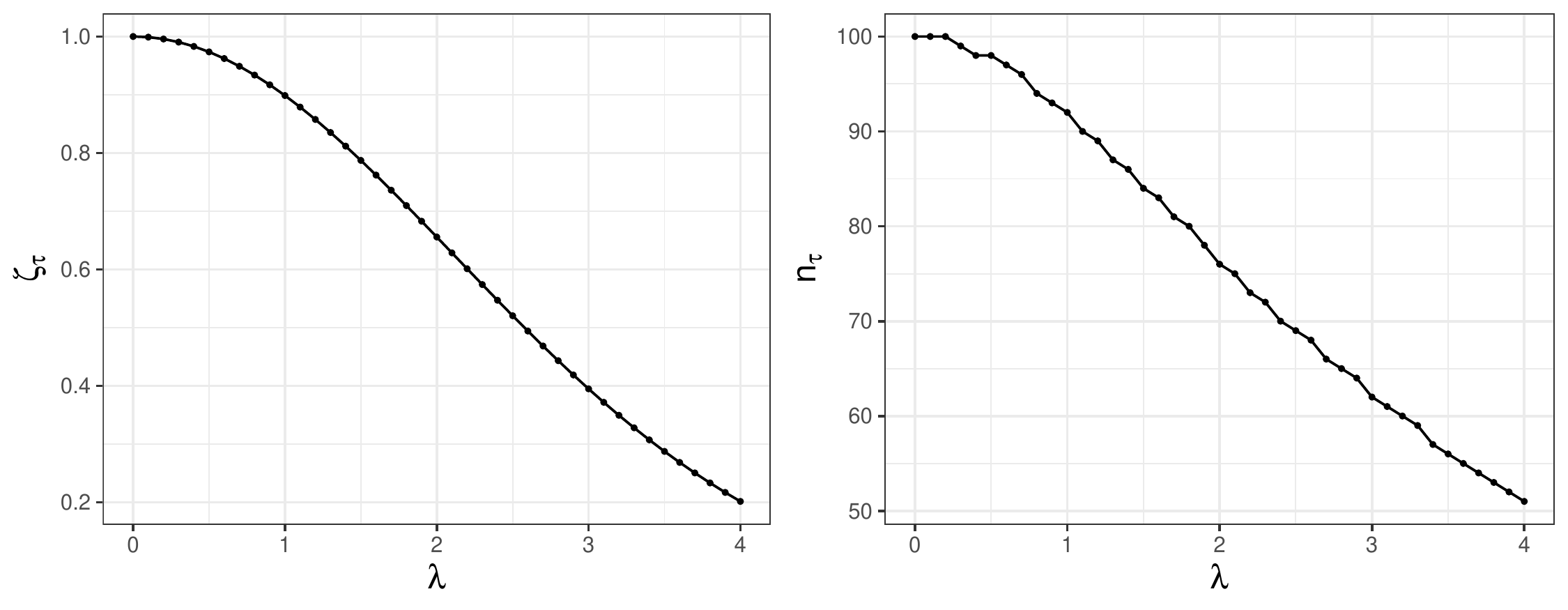}
    \caption{The coefficient $\zeta_\tau$ and effective sample size $n_\tau$ as a function of the skewness parameter $\lambda$ when $n = 100$, $d = 50$, and $\tau(u) = (1 + \vert u \vert)^\lambda$ ($\lambda \geq 0$) with $u \sim \mathcal{N}(0,1)$ and independent of the covariates $\bx$.}
    \label{fig: effective_n}
\end{figure}

In fact, the size of $\zeta_\tau$ and $n_\tau$ are both directly related to the unevenness of $\tau_i$'s. When $w = 1/\tau = 1$, we have $\zeta_\tau = 1$ and $n_\tau = n$. In this case, \eqref{eq: dfR_wls_normal_approx} becomes
\begin{equation} \label{eq: dfR_ols_normal_approx}
    \dfR(p) \approx \frac{p}{2} \left(1 + \frac{n}{n - p - 1}\right),
\end{equation}
which is the approximation of $\dfR$ for subset regression given by \cite{luan2021predictive} under the homoscedastic setting. On the other hand, when $\tau_i$'s are not all equal, $\zeta_\tau < 1$ and $n_\tau \leq n$. In particular, Figure \ref{fig: effective_n} suggests that the more uneven $\tau_i$'s are, the smaller $\zeta_\tau$ and $n_\tau$ will be. Since $n_\tau$ in \eqref{eq: dfR_wls_normal_approx} plays a similar role as $n$ in \eqref{eq: dfR_ols_normal_approx}, we call it the \textit{effective sample size}. When $\tau_i$'s are highly uneven, observations with large variance will receive little weight. The resulting effect is as if the sample size was reduced. 

\subsection{Effect of Weight Schemes} \label{subsec: effect_of_weight}

Since the optimal weight scheme is usually unknown a priori, misspecification of weights could happen in practice. A natural question is then how the model complexity is affected by the choice of weights. In this section, we try to answer this question by fixing the subset of variables for regression and varying the weight scheme $q$ for the weighted least squares model. Since these models are trained under different weight schemes, it is important to assess their complexity using a common criterion. Here, we set $w = 1/\tau$ and define model degrees of freedoms for all weighted least squares models using this weight.

Assume that $\bX \in \R^{n \times p}$ is the design matrix. To avoid confusion, we use $\wErrT_\bX(q)$, $\wErrF_\bX(q)$ and $\wErrR_\bX(q)$ in this section to denote the training and prediction errors of the weighted least squares model $\hat{\bmu}(q)$ with weight scheme $q$. Similarly, we use $\dfF(q)$ and $\dfR(q)$ here to denote the corresponding measures of model degrees of freedom. By replacing $\bW$ with $\bT^{-1}$ in \eqref{eq: dfF_wls_general} and \eqref{eq: dfR_wls_general}, we get $\dfF(q) = \zeta_\tau p$ and 
\begin{equation} \label{eq: dfR_wls_optimal}
    \dfR(q) = \zeta_\tau p + \frac{n}{2} \trace[(\bX^\T \bQ \bX)^{-1} \bX^\T \bQ \bT \bQ \bX (\bX^\T \bQ \bX)^{-1} \bR],
\end{equation}
where $\bR = \bSigma - \bX^\T \bT^{-1} \bX / n$.

For the classical model degrees of freedom $\dfF$, we find that $\dfF$ is constant with respect to the weight scheme $q$. In other words, the average covariance between the fitted values $\hat{\mu}_i$ and the observations $y_i$ is the same for all weighted least squares models. Since $\wErrF_\bX(q) = \wErrT_\bX(q) + (2/n) \bar{w} \sigma_\varepsilon^2 \dfF(q)$, the best weighted least squares model should be the one that minimizes the expected training error $\wErrT_\bX(q)$. Apparently, this minimizer is the model with the optimal weight scheme $q = 1/\tau$, since it is specifically determined to minimize the training error $\wErrT_{\bX, \by}(q) = \Vert \by - \hat{\bmu}(q) \Vert_{\bT^{-1}}^2 / n$ for each $\by$.

For the predictive model degrees of freedom \eqref{eq: dfR_wls_optimal}, the following theorem gives its gradient with respect to $\bq = (q_1,\ldots,q_n)^\T$.
\begin{theorem}\label{eq: differential_dfR_wls_general}
    The gradient of the predictive model degrees of freedom $\dfR$ with respect to the weight vector $\bq = \mathrm{diag}(\bQ)$ is given by
    $$
        \frac{\partial \dfR}{\partial \bq} = n\, \mathrm{diag}\left( (\bI_n - \bH) \bT \bQ \bX (\bX^\T \bQ \bX)^{-1} \bR (\bX^\T \bQ \bX)^{-1} \bX^\T \right),
    $$
    where $\bH = \bX (\bX^\T \bQ \bX)^{-1} \bX^\T \bQ$.
\end{theorem}

For a general weight scheme $q$, it is usually not easy to see directly from $\partial \dfR / \partial \bq$ how $\dfR$ behaves when $q$ is perturbed. However, when $q = 1/\tau$, we have $(\bI_n - \bH) \bT \bQ \bX = (\bI_n - \bH) \bX = \mathbf{O}_{n \times p}$, which further implies that $(\partial \dfR/\partial \bq)\vert_{q = 1/\tau} = \mathbf{0}_{n}$. This indicates that $\dfR$, as a function of the weight scheme $q$, is relatively flat around the optimal weight $1/\tau$. In practice, we usually need to estimate $\tau$ with an estimator $\hat{\tau}$ and use $\hat{q} = 1/\hat{\tau}$ as a weight function for model fitting. A good estimate of $\tau$ guarantees that $\hat{q}$ is close to the optimal weight $1/\tau$. Therefore, it is desirable for $\dfR$ to be flat around the optimal weight in the sense that the discrepancy between the optimal weight scheme $1/\tau$ and a reasonably good estimate of it would make little change in $\dfR$.

\begin{figure}[!t]
    \centering
    \includegraphics[scale = 0.58]{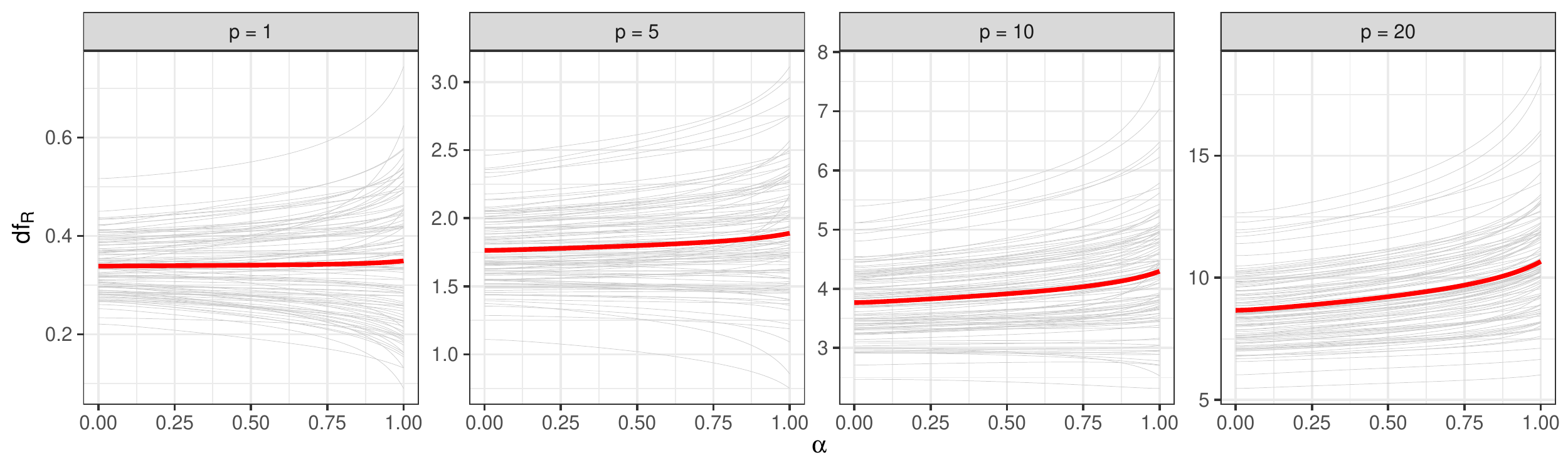}
    \caption{The predictive model degrees of freedom $\dfR$ of the weighted least squares method as a function of the coefficient $\alpha$ when $\bQ = (\alpha \bI_n + (1 - \alpha)\bT)^{-1}$ ($0 \leq \alpha \leq 1$). We assume that $\bx_1,\ldots,\bx_n \sim \mathcal{N}(\mathbf{0},\bI_d)$, $\tau(\bx) \propto (1 + \vert\bx^\T\bbeta\vert)^2$, $\beta_j \propto (1-j/d)^5$ and $\Vert \bbeta \Vert^2 = 10$ with $n = 100$ and $d = 20$. The gray curves represent $\dfR$ based on 100 different realizations of $\bX$, whereas the red curve is the average over these replicates.}
    \label{fig: dfR_vs_alpha}
\end{figure}

A more quantitative study of $\dfR$ is given as follows. Assume that $q$ is of the form
\begin{equation} \label{eq: q_alpha}
    q(\alpha;\bx) = \frac{1}{\alpha + (1 - \alpha) \tau(\bx)}, \text{ for } 0 \leq \alpha \leq 1.
\end{equation}
Then $\E[1/{q(\alpha;\bx)}] = 1 = \E[\tau(\bx)]$ for any $\alpha$. Note that $q = 1/\tau$ when $\alpha = 0$ and $q = 1$ when $\alpha = 1$. Thus, \eqref{eq: q_alpha} constructs a path from the optimal to the equal weight scheme. To examine $\dfR$ as a function of $\alpha$ more closely, we consider a simulation study below. We set $n = 100$ and $d = 20$, and assume that $\bx_1,\ldots,\bx_n \sim \mathcal{N}(\mathbf{0},\bI_d)$. Also assume that $\tau(\bx) \propto (1 + \vert \bx^\top\bbeta \vert)^2$ with $\beta_j \propto (1-j/d)^5$ and $\Vert \bbeta \Vert^2 = 10$. We generate 100 replicates of $\bX$. Figure \ref{fig: dfR_vs_alpha} demonstrates the predictive model degrees of freedom of models with the first $p$ variables for $p = 1,\,5,\,10$ and $20$. We can see that $\dfR$ can be either increasing or decreasing in $\alpha$. However, as the number of variables increases, $\dfR$ is more likely to be an increasing function of $\alpha$. Note that, as more relevant variables are added, the weighted prediction errors $(y_\ast - \hat{\mu}_\ast)/\sqrt{\tau_\ast}$ of the model with the optimal weight scheme are generally homoscedastic, while those of the other models are not. Consequently, the variance of the predicted values produced by the model with optimal weights is the lowest on average. Since the predictive model degrees of freedom is defined through the excess variance, the optimal weight scheme generally leads to the simplest model.

On the other hand, though the weight scheme $q$ for model fitting affects the predictive model degrees of freedom, our results in Figure \ref{fig: dfR_vs_alpha} suggest that its impact may not be big enough to produce appreciable range of difference. When the model size is not too large, which is usually the case in many real applications, the difference in $\dfR$ due to the weight scheme $q$ for model fitting is fairly small or even trivial, as we can see from the figure. In addition, the expected optimism \eqref{eq: OptR_weighted_decomposition} implies that
\begin{equation} \label{eq: ErrR_decomposition}
    \wErrR_\bX = \wErrT_\bX + \sigma_\varepsilon^2 \left(\E(w_\ast \tau_\ast) - \frac{1}{n}\sum_{i=1}^n w_i \tau_i \right) + \Delta B_\bX + \frac{2}{n} \bar{w} \sigma_\varepsilon^2 \dfR.
\end{equation}
Since $\dfR$ contributes to the prediction error $\wErrR_\bX$ through $(2/n)\bar{w} \sigma_\varepsilon^2 \dfR$, the impact of the difference in $\dfR$ for different weighted least squares models will be reduced when the sample size $n$ is large. By contrast, the choice of $q$ may have a much larger impact on the training error $\wErrT_\bX$ or excess bias $\Delta B_\bX$. In Section \ref{sec: pred_err_est}, we will take these terms into account and provide a risk estimator for model evaluation and comparison.

\section{Prediction Error Estimation} \label{sec: pred_err_est}

In this section, we aim to construct an estimator for the out-of-sample prediction error based on \eqref{eq: ErrR_decomposition}. For the excess irreducible error, the law of large numbers guarantees that $\E(w_\ast \tau_\ast) - \sum_{i = 1}^n w_i \tau_i / n$ is generally negligible when $n$ is large. The excess bias $\Delta B_\bX$ depends on the unknown true mean function $\mu(\cdot;\bbeta)$, and thus needs to be estimated. The predictive model degrees of freedom depends on the distribution of $\bx_\ast$ through $\E(w_\ast \Vert \bh_\ast \Vert_\bT^2 \vert \bX)$. In practice, this expectation also needs to be estimated.

For simplicity, we assume that $\tau$ is either known a priori or can be estimated by $\hat{\tau}$ using some independent data (e.g., historical data) throughout this section. Further, we assume that $q = w \propto 1/\tau$ (or $1/\hat{\tau}$). 

\subsection{Estimation of \texorpdfstring{$\Delta B_\bX$}{\textDelta B}}
We first consider estimating the excess bias $\Delta B_\bX$. Here, we adjust with weights the method proposed in \cite{luan2021predictive} that uses the leave-one-out cross validation (LOOCV) trick. Let $\bX^{-i}$, $\by^{-i}$ and $\bmu^{-i}$ be the corresponding terms with the $i$th case deleted. Write $\bH = (\bh_1,\ldots,\bh_n)^\T$ and define $\bh^{-i}_i$ to be the hat vector at $\bx_{i}$ based on $\bX^{-i}$. Conceptually, when $n$ is large, we have
$$
	\E[w_\ast (\mu_\ast - \bh_\ast^\T\bmu)^2 \vert \bX] \approx \frac{1}{n} \sum_{i=1}^n w_i (\mu_i - (\bh^{-i}_i)^\T \bmu^{-i})^2.
$$
For linear procedures, the right hand side of the above equation can be evaluated efficiently due to the LOOCV identity $y_i - (\bh_i^{-i})^\T \by^{-i} = (y_i - \bh_i^\T \by) / (1 - h_{ii})$ \citep{craven1978smoothing}. Since $\bh_i$ doesn't depend on either $\bmu$ or $\by$, the identity still holds if we replace $y_i$ with $\mu_i$ and $\by^{-i}$ with $\bmu^{-i}$ respectively. Then, we have
$$
\begin{aligned}
	\Delta B_\bX
	& \approx \frac{1}{n}\sum_{i=1}^n w_i [(\mu_i - (\bh^{-i}_i)^\T \bmu^{-i})^2 - (\mu_i - \bh_i^\T \bmu)^2]\\
	& = \frac{1}{n}\sum_{i=1}^n w_i \left[\frac{(\mu_i - \bh_i^\T\bmu)^2}{(1-h_{ii})^2} - (\mu_i - \bh_i^\T\bmu)^2\right]\\
	& = \frac{1}{n}\bmu^\T \mathbf{A} \bmu \\
	& = \frac{1}{n}\by^\T \mathbf{A} \by - \frac{2}{n}\bmu^\T \mathbf{A} \bm{\varepsilon} - \frac{1}{n}\bm{\varepsilon}^\T \mathbf{A} \bm{\varepsilon},
\end{aligned}
$$
where $\mathbf{A}= (\bI_n - \bH)^\T\mathbf{D}(\bI_n - \bH)$ and $\mathbf{D} = \mathrm{diag}\left(w_i/(1-h_{ii})^2 - w_i\right)$. Since $\E(\bmu^\T \mathbf{A}\bm{\varepsilon}\vert \bX) = 0$ and $\E(\bm{\varepsilon}^\T \mathbf{A} \bm{\varepsilon}\vert \bX) = \sigma_\varepsilon^2 \trace(\mathbf{A}\bT)$, we can estimate $\Delta B_\bX$ by
$$
	\hat{\delta} = \frac{1}{n} \by^\T \mathbf{A} \by - \frac{1}{n}\sigma_\varepsilon^2 \trace(\mathbf{A}\bT).
$$
Assume $\dfR$ is known.  Then \eqref{eq: ErrR_decomposition} implies that an estimator of $\wErrR_\bX$ is given by
\begin{equation} \label{eq: ErrR_hat}
    \widehat{\wErrR} = \wErrT_{\bX,\by} + \hat{\delta} + \frac{2}{n} \bar{w} \sigma_\varepsilon^2 \dfR = \widehat{\wErrR}_{\rm loocv} + \frac{1}{n} \sigma_\varepsilon^2 \xi_\bX,
\end{equation}
where $\widehat{\wErrR}_{\rm loocv} = \sum_{i=1}^n w_i (y_i - (\bh_i^{-i})^\T \by^{-i})^2 / n$ is the LOOCV error and $\xi_\bX = 2\bar{w}\dfR - \trace(\mathbf{A} \bT)$. Thus, $\widehat{\wErrR}$ is an adjusted version of the LOOCV error. In fact, this adjustment has a close connection with the variance of prediction errors in the full model and leave-one-out models, which we summarize in the following theorem. 
\begin{theorem} \label{thm: adjustment_interpretation}
    The adjustment term $\frac{1}{n}\sigma_\varepsilon^2 \xi_\bX$ in \eqref{eq: ErrR_hat} can be expressed as
    $$
        \frac{1}{n} \sigma_\varepsilon^2 \xi_\bX = \E[w_\ast \var(\hat{\varepsilon}_\ast \vert \bx_\ast, \bX) \vert \bX] - \frac{1}{n} \sum_{i=1}^n w_i \var(\hat{\varepsilon}^{-i}_i \vert \bX) - \sigma_\varepsilon^2 \left( \E(w_\ast \tau_\ast) - \frac{1}{n}\sum_{i=1}^n w_i \tau_i\right),
    $$
    where $\hat{\varepsilon}_\ast$ is the prediction error of the full model at $\bx_\ast$, and $\hat{\varepsilon}_i^{-i}$ is that of the leave-the-$i$th-out model at $\bx_i$.
\end{theorem}

Note that the third term on the right hand side of the expression above is the excess irreducible error, which is approximately 0 when $n$ is large. Therefore, $\frac{1}{n}\sigma_\varepsilon^2 \xi_\bx $ can be approximated by
$$
    \frac{1}{n} \sigma_\varepsilon^2 \xi_\bX \approx \E[w_\ast \var(\hat{\varepsilon}_\ast \vert \bx_\ast, \bX) \vert \bX] - \frac{1}{n} \sum_{i=1}^n w_i \var(\hat{\varepsilon}^{-i}_i \vert \bX).
$$
This indicates that the adjustment term $\frac{1}{n}\sigma_\varepsilon^2 \xi_\bX$ is a measure of the difference between the full model and leave-one-out models in the variance of prediction errors weighted by the case weights $w_1,\ldots,w_n$. In addition, as shown in \cite{luan2021predictive}, the adjustment term is generally negligible when $n$ is large and $p/n$ is relatively small. But when $p$ is close to $n$, $\trace(\bA\bT)$ could be very large and variable due to the leverages $h_{ii}$ approaching 1. This may make the adjustment negative and large in magnitude, and further cause $\widehat{\wErrR}$ to be negative. Since the term $\trace(\bA \bT)$ in $\xi_\bX$ comes from the excess bias estimator $\hat{\delta}$, to circumvent this issue, they suggested to use $\hat{\delta}_+ = \max(\hat{\delta}, 0)$ in place of $\hat{\delta}$ in \eqref{eq: ErrR_hat}. We will follow this adjustment in the rest of the paper.

\subsection{Estimation of \texorpdfstring{$\dfR$}{dfR}} \label{subsec: estimate_dfR_summary}
In practice, the distribution of $\bx_\ast$ is usually unknown, which makes $\E(w_\ast \Vert \bh_\ast \Vert_\bT^2 \vert \bX)$ and $\dfR$ not directly available. To estimate $\E(w_\ast \Vert \bh_\ast \Vert_\bT^2 \vert \bX)$, one can again use the LOOCV trick and approximate the expectation by $\sum_{i=1}^n w_i \Vert \bh_i^{-i} \Vert_\bT^2 / n$. Unfortunately, unlike estimating the excess bias $\Delta B_\bX$, this approach inevitably requires fitting $n$ leave-one-out models, which could be highly inefficient when $n$ is large.

Another approach is to estimate the expectation $\E(w_\ast \Vert \bh_\ast \Vert_\bT^2 \vert \bX)$ using some synthetic data generated based on the training data $\bX$ that might mimic genuine random samples from the joint distribution of $\bx$. Assume that $B$ synthetic examples are generated, denoted by $\bx_\ast^{(1)}, \ldots, \bx_\ast^{(B)}$. Then we can estimate the predictive model degrees of freedom by
$$
    \widehat{\mathrm{df}}_\mathrm{R} = \dfF + \frac{n}{2\bar{w}}\left(\frac{1}{B} \sum_{b=1}^B w_\ast^{(b)} \Vert \bh_\ast^{(b)} \Vert_\bT^2 - \frac{1}{n}\trace(\bW \bH \bT \bH^\T)\right),
$$
where $w_\ast^{(b)} = w(\bx_\ast^{(b)})$ and $\bh_\ast^{(b)}$ is the hat vector at $\bx_\ast^{(b)}$ ($b = 1,\ldots,B$).

In this work, we consider two different procedures to generate synthetic data: the Classification and Regression Trees (CART) data synthesis method \citep{reiter2005using} and the Naive Bayes Estimation (NBE) learning algorithm \citep{lowd2005naive}. One big advantage of these two methods is that they are efficient and readily available without any extra model tuning. For a given set of variables $\bx = (x_1,\ldots, x_d)$, the former seeks to estimate the full conditional distributions $p(x_j \vert x_1,\ldots, x_{j-1}, x_{j+1}, x_d)$ using CART \citep{breiman2017classification} and sample sequentially from these distributions. The NBE algorithm, on the other hand, takes a parametric approach by assuming conditional independence among $X_j$'s given a latent variable $Z$ and estimating the conditional joint distribution $p(\bx\vert z)$ by a mixture of models (Gaussian for continuous variables and multinomial for categorical ones) using the EM algorithm \citep{dempster1977maximum}. In Section \ref{sec: numerical_study}, we compare the performance of these two methods in estimation of $\dfR$ through some simulation studies.

\begin{figure}[!t]
    \centering
    \includegraphics[scale = 0.6]{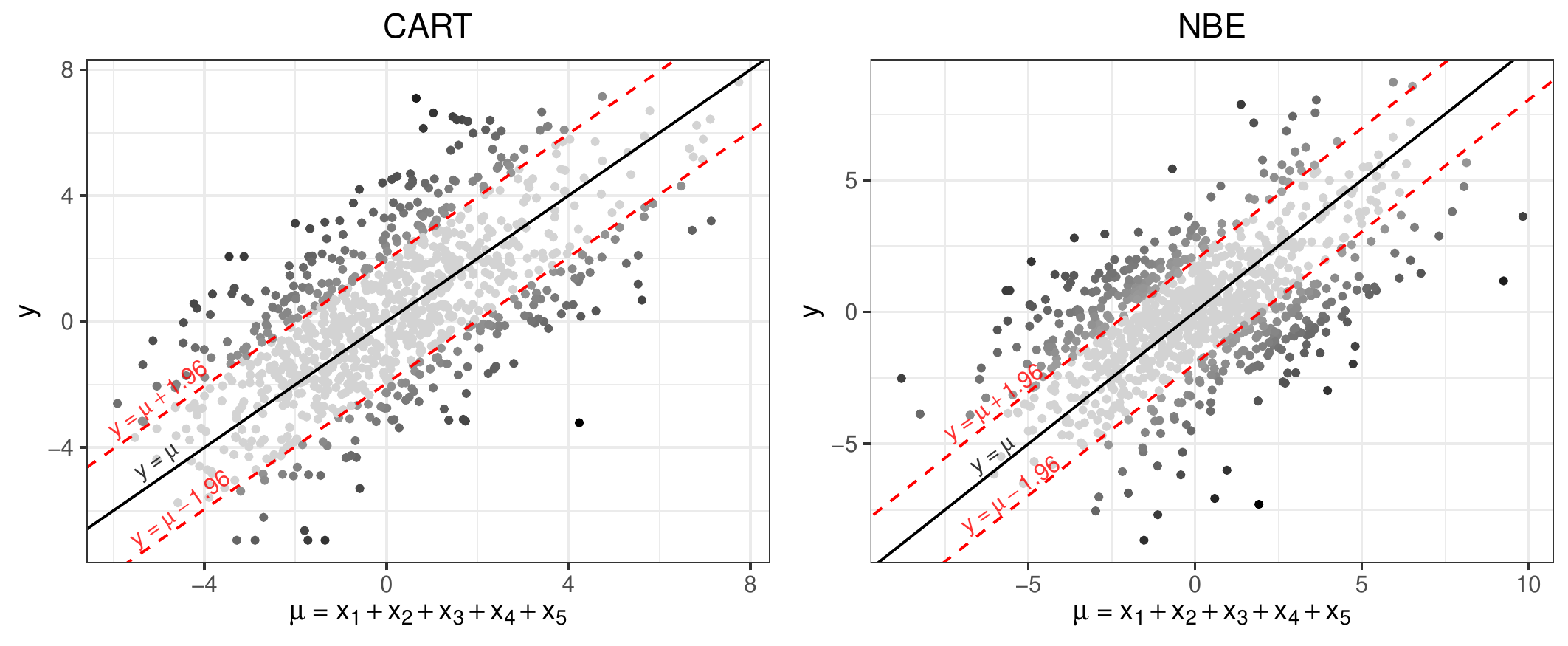}
    \caption{Response $y$ against its mean $\mu = \mu(\bx)$ from synthetic data generated based on a fixed training dataset of size 200. The true data generating mechanism is $y_i = \sum_{j=1}^5 x_{ij} + \varepsilon_i$, where $x_{ij}$'s and $\varepsilon_i$ are i.i.d.\ $\mathcal{N}(0,1)$. 1,000 synthetic examples are generated for each method.}
    \label{fig: syn_simXY}
\end{figure}

Here, one question remains to be answered: Can we synthesize $\lbrace (\bx_\ast, y_\ast)\rbrace_{b = 1,\ldots,B}$ directly and evaluate out-of-sample prediction error based on the synthetic data? To study its feasibility, we consider a simple example below. Assume $y_i = \sum_{j=1}^5 x_{ij} + \varepsilon_i$, where $x_{ij}$'s and $\varepsilon_i$ are i.i.d.\ $\mathcal{N}(0,1)$. We set $n = 200$ and $B = \text{1,000}$, and generate 500 synthetic datasets based on a fixed training dataset $(\bX,\by)$ using both the NBE and CART algorithms. We apply the ordinary least squares method with all variables included to fit a linear model. It is easy to show that the true out-of-sample prediction error is given by $\ErrR_\bX = 1 + \trace[(\bX^\T \bX)^{-1}] \approx 1.027$. However, we find that both data synthesizing methods overestimate the true prediction error, with an average of 3.872 (SD = 0.185) for the CART method and 5.878  (SD = 0.270) for the NBE algorithm. For each method, a scatterplot of the synthesized $y_\ast$'s against $x_\ast$'s is given in Figure \ref{fig: syn_simXY}. While the linear pattern between $x_\ast$ and $y_\ast$ still remains as expected, $y_\ast$ clearly has variance greater than 1. This explains the inflated prediction error estimates. Our numerical results indicate that such over-dispersion generally becomes more notable as the signal-to-noise ratio increases. Structurally, the full conditional distributions estimated using CART and the joint distribution estimated through the NBE algorithm under the assumption of conditional independence may not provide a good approximation of the true joint distribution of $(\bx,y)$. Therefore, it is not recommended to estimate the prediction error based on synthetic data directly. 

On the other hand, this result doesn't invalidate our approach to estimating $\dfR$ using synthetic data. In fact, our numerical results in Section \ref{subsec: simulations} demonstrate great robustness of $\widehat{\mathrm{df}}_{\rm R}$ to synthetic data even when all features are moderately correlated. Also, note that the estimated predictive model degrees of freedom $\widehat{\mathrm{df}}_{\rm R}$ contributes to the prediction error estimate $\widehat{\wErrR}$ through the term $(2/n) \bar{w} \sigma_\varepsilon^2 \widehat{\mathrm{df}}_{\rm R}$. When $n$ is large, the factor $2/n$ in front can effectively alleviate the impact of discrepancy between the true and estimated joint covariate distribution. Besides, we can further reduce this impact by wisely selecting variables for data synthesis. For example, when multiple functions of a certain variable $X_j$ are used as predictors in a model (e.g., $X_j,\ldots, X_j^k$ in polynomial regression), it is good practice to synthesize $X_j$ first and propagate it to other transformations.

\section{Numerical Studies} \label{sec: numerical_study}

In this section, we evaluate the performance of the estimators developed in Section \ref{sec: pred_err_est} through numerical studies. 

\subsection{Simulations} \label{subsec: simulations}

Throughout our experiments, we set $n = 60$ and generate 20 continuous variables and 20 categorical variables in total ($d = 40$). Let $\mathbf{V} = (V_1,\ldots,V_{20})^\T \in \R^{20}$ and $\mathbf{C} = (C_1,\ldots,C_{20})^\T \in \lbrace -1,+1\rbrace^{20}$ be the sets of continuous and binary categorical variables respectively. To specify relationship among features with different level of complexity, we consider the following two settings. 
\begin{itemize}
    \item Single component. Assume that $\bV \sim \mathcal{N}(\mathbf{0}, \bSigma_\bV)$ is independent of $\mathbf{C}$, and $C_1,\ldots,C_{20}$ are i.i.d.\ with $\prob(C_1 = +1) = \prob(C_1 = -1) = 1/2$.
    \item Mixture of components. Let $Z \in \lbrace -1,+1\rbrace$ be a latent random variable such that $\prob(Z = +1) = \prob(Z = -1) = 1/2$. Assume that $\bV \vert Z \sim \mathcal{N}(Z(\log 1,\ldots, \log 20)^\T, \bSigma_\bV)$, $\prob(C_1 = +1) = 1/2$ and $\prob(C_j = +1 \vert C_{j-1}, Z) = 1/2 + Z C_{j-1}/(4 + \log j)$ for $j = 2,\ldots, 20$.
\end{itemize}
For the variance-covariance matrix $\bSigma_\bV$, we consider three different cases: 1) I.I.D.: $\bSigma_\bV = \bI_{20}$; 2) Equicorrelation: $\bSigma_\bV = \frac{1}{2} \bI_{20} + \frac{1}{2}\mathbf{1}_{20}\mathbf{1}_{20}^\T$; and 3) Autoregressive correlation: $\bSigma_\bV = (2^{-\vert i-j \vert})_{i,j=1,\ldots,20}$. Let $\bX = (\bx_1,\ldots,\bx_n)^\T$ be the design matrix. We assume that the true data generating mechanism $\mu(\cdot;\bbeta)$ is linear, i.e., $\mu(\bx_i; \bbeta) = \sum_{j=1}^d \beta_j x_{ij}$. For convenience, we use $\beta^{\rm v}_j$ and $\beta^{\rm c}_j$ to denote the coefficient $\beta_j$'s for the corresponding continuous and categorical variables respectively. We set $\beta^{\rm v}_j = \beta^{\rm c}_j/2 = a \left(1 - j/20\right)^5$ and choose $a$ such that $\sum_{j = 1}^{20} (\beta^{\rm v}_j)^2 = 5$. In subset regression, we add continuous and categorical variables alternatingly in descending order of their coefficients, i.e., in the order $V_1, C_1, \ldots,V_{20}, C_{20}$.

For the error variances, we assume $\sigma_\varepsilon^2 = 1$ and $\tau(\bx) \propto (1 + \vert \mu(\bx;\bbeta)\vert)^2$. We assume $\tau(\bx)$ can be obtained for each $\bx$ without knowing the mean function $\mu$. We use $w = 1/\tau$ as the weight function when fitting weighted least squares models.

\subsubsection{Estimation of \texorpdfstring{$\dfR$}{dfR}}
\begin{figure}[t!]
	\centering
	\includegraphics[scale = 0.5]{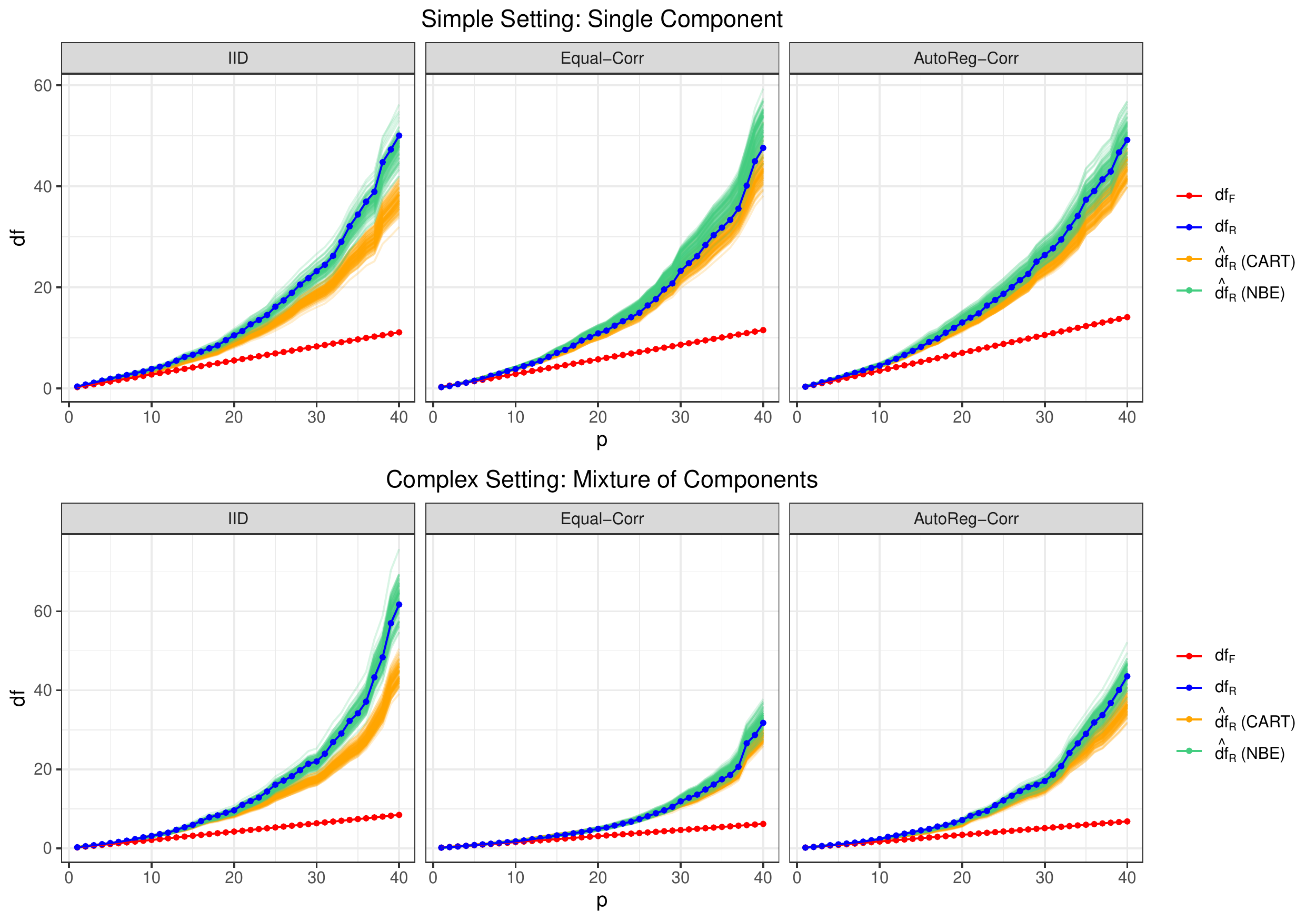}
	\caption{Estimated predicted model degrees of freedom $\widehat{\mathrm{df}}_{\rm R}$ based on synthetic data generated using CART and NBE algorithm. For each method, 100 synthetic datasets of size 1,000 are generated based on a single replicate of $\bX$ of size $n = 60$ and $d = 40$ (20 continuous and 20 categorical). The true mean function is linear in $\bx$, and continuous and categorical variables are added alternatingly in the descending order of their coefficients.}
	\label{fig: df_wls_est}
\end{figure}

We consider the CART method and the NBE learning algorithm described in Section \ref{subsec: estimate_dfR_summary} to generate synthetic data. Figure \ref{fig: df_wls_est} illustrates and compares the results using these two procedures based on a random realization of $\bX$. 
Each method is applied 100 times to generate different synthetic datasets used for estimation of $\dfR$. We can see a consistent pattern across all 6 examined cases. Both methods exhibit stability with respect to the choice of synthetic data, though the variance slightly inflates as the subset size $p$ increases. The CART method tends to underestimate the true predictive model degrees of freedom, especially under the I.I.D.\ setting. These estimates lie between the predictive model degrees of freedom $\dfR$ and the classical model degrees of freedom $\dfF$. 
This suggests that the synthetic data generated using this method do not fully mimic out-of-sample examples but act like a ``hybrid'' of out-of-sample and in-sample data. This is not surprising since the method samples from a set of full conditional distributions as mentioned in Section \ref{subsec: estimate_dfR_summary}, which does not guarantee that the resulting synthetic data follow the same joint distribution as the training data. On the other hand, despite a slight tendency of overestimation in the equicorrelation setting, the NBE algorithm performs remarkably well in estimating the true predictive model degrees of freedom, even when the data do not meet the ``naive'' assumption that all features are conditionally independent given the latent cluster. \cite{lowd2005naive} pointed out that the naive Bayes mixture model is a special case of Bayesian networks, which are commonly used in general high-dimensional probability representation and estimation. In addition, the naive Bayes mixture model is readily available and orders of magnitude faster than Bayesian networks for sampling and inference. This enables us to estimate the predictive model degrees of freedom in an accurate and efficient fashion.

\subsubsection{Estimation of \texorpdfstring{$\Delta B_\bX$}{\textDelta B} and \texorpdfstring{$\wErrR_\bX$}{wErrR}} \label{subsubsec: estimate_ErrR}

For simplicity, we only consider the setting of a mixture of components with i.i.d.\ continuous covariates ($\bSigma_\bV = \bI_{20}$). We generate 500 replicates of $(\bX, \by)$. Based on each replicate, the estimators $\hat{\delta}$ and $\widehat{\wErrR}$ are calculated along with the true excess bias $\Delta B_\bX$ and prediction error $\wErrR_\bX$. In general, the excess bias estimator is unbiased when $p$ is small, but becomes slightly biased upward as $p$ increases due to the truncation of $\hat{\delta}$, $\hat{\delta}_+ = \max(\hat{\delta},0)$ when $\hat{\delta} < 0$. 
\begin{figure}[!t]
    \centering
	\includegraphics[scale = 0.57]{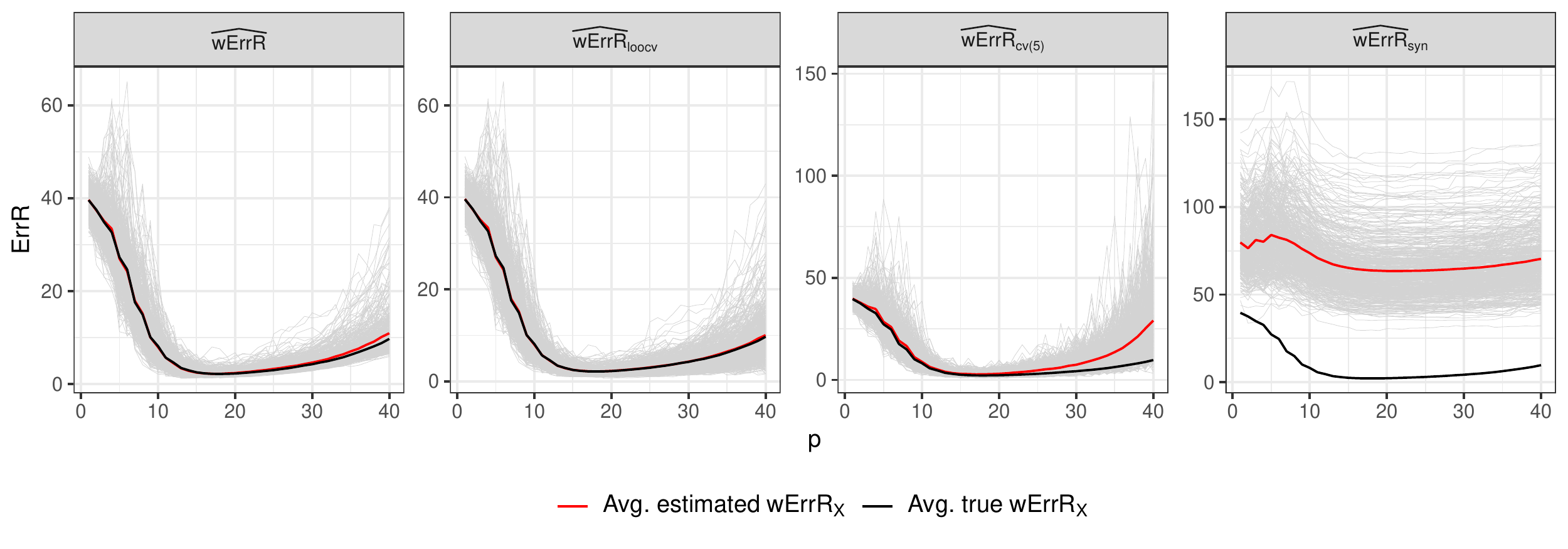}
	\caption{Comparison of the four estimators $\widehat{\wErrR}$, $\widehat{\wErrR}_{\rm loocv}$, $\widehat{\wErrR}_{\rm cv(5)}$ and $\widehat{\wErrR}_{\rm syn}$ in prediction error estimation. The gray lines are the risk estimates based on each of the 500 random replicates of $(\bX, \by)$ generated under the mixture components setting with $\bSigma_\bV = \bI_{20}$. The red line marks the average of these estimates, whereas the black line shows the average of the true prediction ever over all replicates.}
	\label{fig: ErrR_wls_est}
\end{figure}

With the corrected excess bias estimator $\hat{\delta}_+$ and the estimated predictive model degrees $\widehat{\mathrm{df}}_{\rm R}$ based on the NBE algorithm, we form the prediction error estimator $\widehat{\wErrR}$. For comparison, we also include three other estimators: the leave-one-out cross validation error $\widehat{\wErrR}_{\rm loocv}$, the 5-fold cross validation error $\widehat{\wErrR}_{\rm cv(5)}$, and the prediction error estimator based on synthetic data, denoted by $\widehat{\wErrR}_{\rm syn}$. The results are illustrated in Figure \ref{fig: ErrR_wls_est}. We can see that our estimator is very similar to the LOOCV error when $p$ is relatively small, but has much lower variance when $p$ gets large. This results from the use of $\hat{\delta}_+$ that trades off a small increase in bias for a large reduction in variance. Note that 5-fold cross validation doesn't work very well in our example, since the training folds contain only 48 observations but a maximum of 40 covariates, which inevitably inflates bias and variance simultaneously. $\widehat{\wErrR}_{\rm syn}$ performs the worst due to the over-dispersion issue in the synthetic data.
\begin{figure}[!t]
    \centering
	\includegraphics[scale = 0.53]{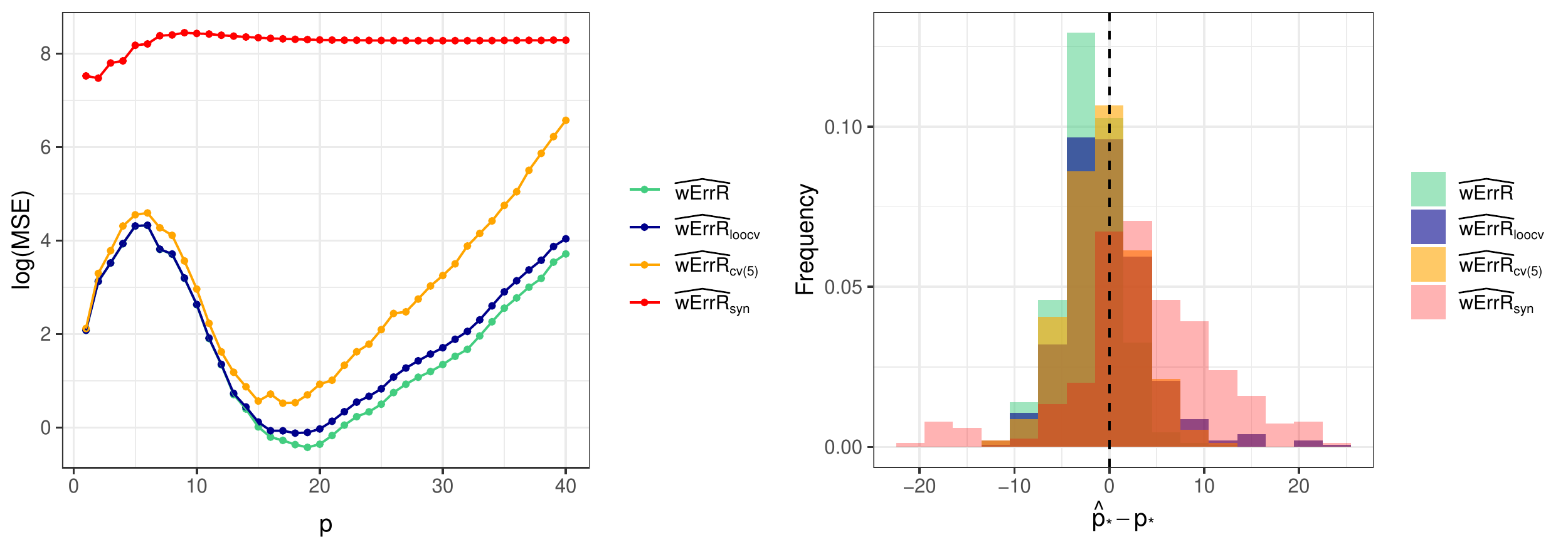}
	\caption{Comparison of the four prediction error estimators $\widehat{\wErrR}$, $\widehat{\wErrR}_{\rm loocv}$, $\widehat{\wErrR}_{\rm cv(5)}$ and $\widehat{\wErrR}_{\rm syn}$ in terms of mean squared error (left) and model selection (right) based on 500 replicates of $(\bX, \by)$ generated under the mixture components setting with $\bSigma_\bV = \bI_{20}$. For each $p$, the mean squared error of an estimator is the average squared difference between the estimate and true prediction error over all replicates. For each replicate, we define $p_\ast$ and $\hat{p}_\ast$, respectively, as the minimizer of the true and estimated prediction errors.}
	\label{fig: ErrR_wls_est_compr}
\end{figure}

A more quantitative comparison of the four prediction error estimators is given in Figure \ref{fig: ErrR_wls_est_compr}. The left panel shows the mean squared error (on log scale) of each estimator as a function of $p$. Our estimator $\widehat{\wErrR}$ beats the other three due to its much smaller variance. In terms of model selection, let $p_\ast$ and $\hat{p}_\ast$ be the minimizers of the true and estimated prediction errors respectively for each replicate. The right panel shows the histogram of $\hat{p}_\ast - p_\ast$. We can see that $\widehat{\wErrR}$ spreads the least, though its center is a bit off to the left of 0, meaning that it tends to select a slightly simpler model than the true model. Since out-of-sample prediction generally involves more uncertainty than in-sample prediction, it is reasonable to be slightly conservative in model selection. On the other hand, the other three estimators are more likely to select larger models.
\begin{figure}[!t]
    \centering
	\includegraphics[scale = 0.5]{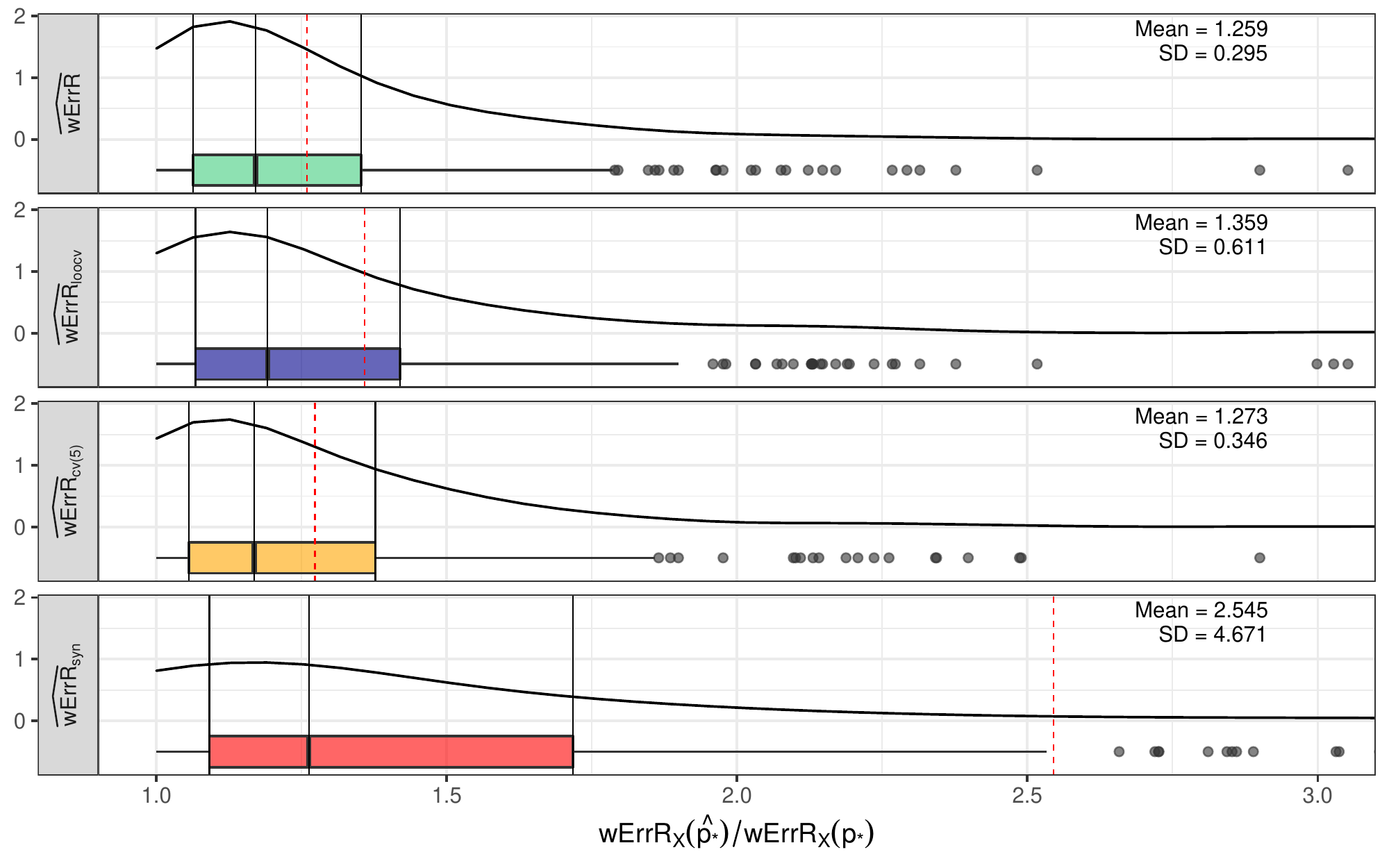}
	\caption{Comparison of models selected by the four prediction error estimators $\widehat{\wErrR}$, $\widehat{\wErrR}_{\rm loocv}$, $\widehat{\wErrR}_{\rm cv(5)}$ and $\widehat{\wErrR}_{\rm syn}$ based on 500 replicates of $(\bX, \by)$ generated under the mixture components setting with $\bSigma_\bV = \bI_{20}$. The red dashed lines show the average ratio $\wErrR_\bX(\hat{p}_\ast)/\wErrR_\bX(p_\ast)$.}
	\label{fig: ErrR_wls_est_compr2}
\end{figure}

A simpler model is desirable in the sense that it is usually computationally more efficient and has better interpretability. But if prediction is a top priority, how much predictive accuracy do we lose in pursuit of such model parsimony? To answer this question, we look at the ratio of the true prediction error of the selected model to that of the optimal one, defined as $\wErrR_\bX (\hat{p}_\ast)/\wErrR_\bX (p_\ast)$. For a good model selection criterion, this ratio is expected to be close to 1. As illustrated in Figure \ref{fig: ErrR_wls_est_compr2}, our risk estimator $\widehat{\wErrR}$ exhibits the best performance in both the mean and standard deviation of the ratio. The 5-fold cross-validation error $\widehat{\wErrR}_{\rm cv(5)}$ also works comparably well, but its variance is slightly larger. By contrast, the risk estimator based on synthetic data is very unstable and can easily result in poor models.

\subsubsection{WLS vs. OLS} \label{subsubsec: wls_vs_ols}
So far in our numerical studies, we have assumed that the the weight function is correctly specified. In practice, however, since error variances are usually unknown, case weights need to be estimated and could be misspecified, which may consequently affect statistical inference and decision. In this section, we examine the effect of weight misspecification for the least squares method. In particular, we focus on model selection and compare two extreme weight schemes: the optimal weight and equal weight schemes under both in-sample and out-of-sample prediction settings. For simplicity, we still assume that the data are generated under the mixture components setting with $\bSigma_\bV = \bI_{20}$, and that $\tau(\bx) \propto (1 + \vert \mu(\bx;\bbeta)\vert)^\lambda$ with $\lambda = 2$. 

\begin{figure}[t!]
    \centering
    \includegraphics[scale = 0.6]{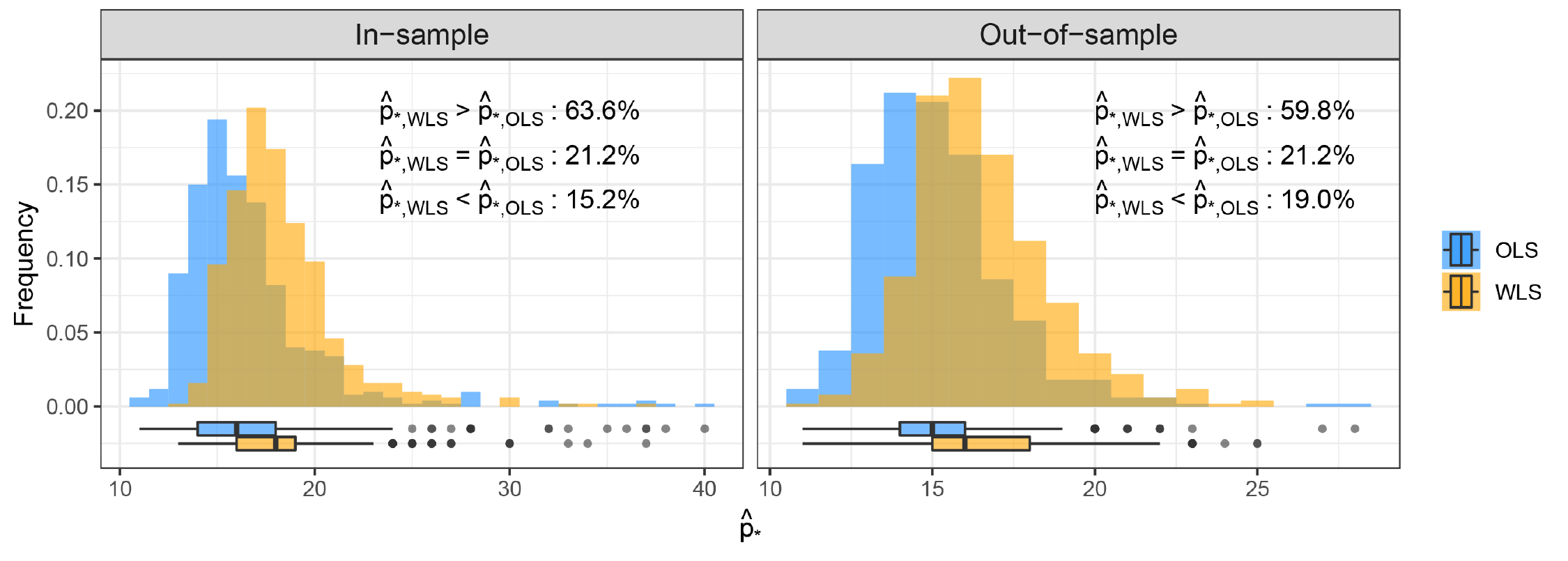}
    \caption{Comparison of the selected model size $\hat{p}_\ast$ between the ordinary least squares model and the weighted least squares model with the optimal weight scheme based on 500 replicates of $(\bX, \by)$ generated under the mixture components setting with $\bSigma_\bV = \bI_{20}$ and $\tau(\bx) \propto (1 + \vert \mu(\bx;\bbeta)\vert)^2$. $\widehat{\wErrF}$ and $\widehat{\wErrR}$ are used as the model selection criterion in the in-sample and out-of-sample prediction settings respectively. }
    \label{fig: model_selection_ols_vs_wls}
\end{figure}

We use $\widehat{\wErrF}$ and $\widehat{\wErrR}$ as a model selection criterion in the in-sample and out-of-sample prediction settings respectively. As shown in Figure \ref{fig: model_selection_ols_vs_wls}, the ordinary least squares model tends to select fewer variables than the weighted least squares model with the optimal weight scheme in both scenarios. This difference in the size of the selected model is generally the same for in-sample and out-of-sample predictions. In terms of the variance of the estimated model size $\hat{p}_\ast$, the two methods are quite close, with the weighted least squares model being slightly more variable under the out-of-sample prediction setting. 

Since $\widehat{\wErrR}$ tends to favor a slightly more parsimonious model as we have shown in Section \ref{subsubsec: estimate_ErrR}, misspecification of weights could lead to an oversimplified model. This could further result in a nontrivial increase in the prediction error if the true risk curve to the left of the optimal model size is steep. In our example, the in-sample and out-of-sample prediction errors could be inflated by 12\% and 50\% respectively on average due to weight misspecification. This issue will become more substantial with more uneven error variances.

\subsection{Real Data Analysis} \label{subsec: rda}
We consider two real data problems in this section. 

\subsubsection{Cancer Mortality Data} \label{subsubsec: rda_cancer}

\begin{figure}[!t]
	\centering
	\includegraphics[scale = 0.48]{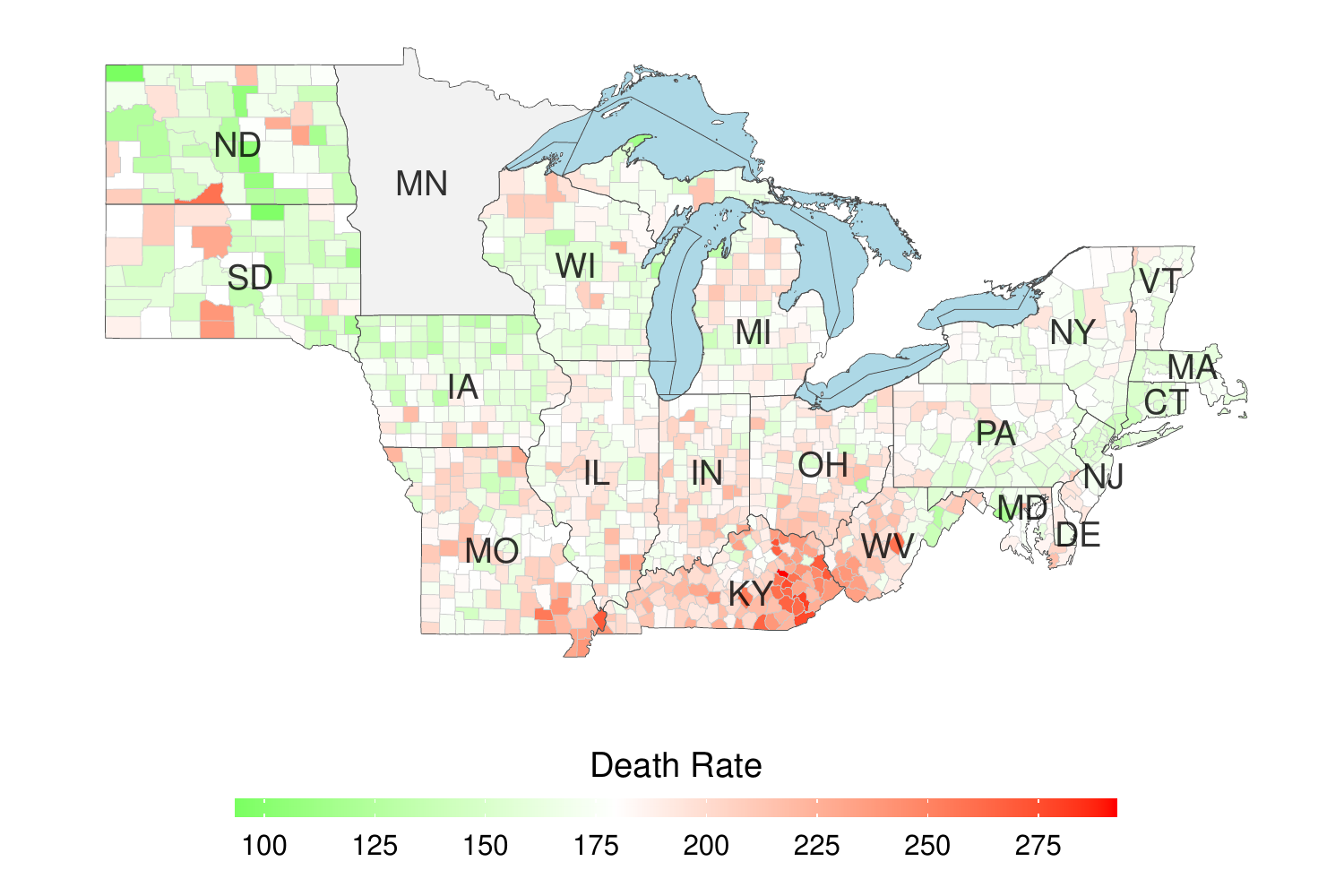}
	\includegraphics[scale = 0.5]{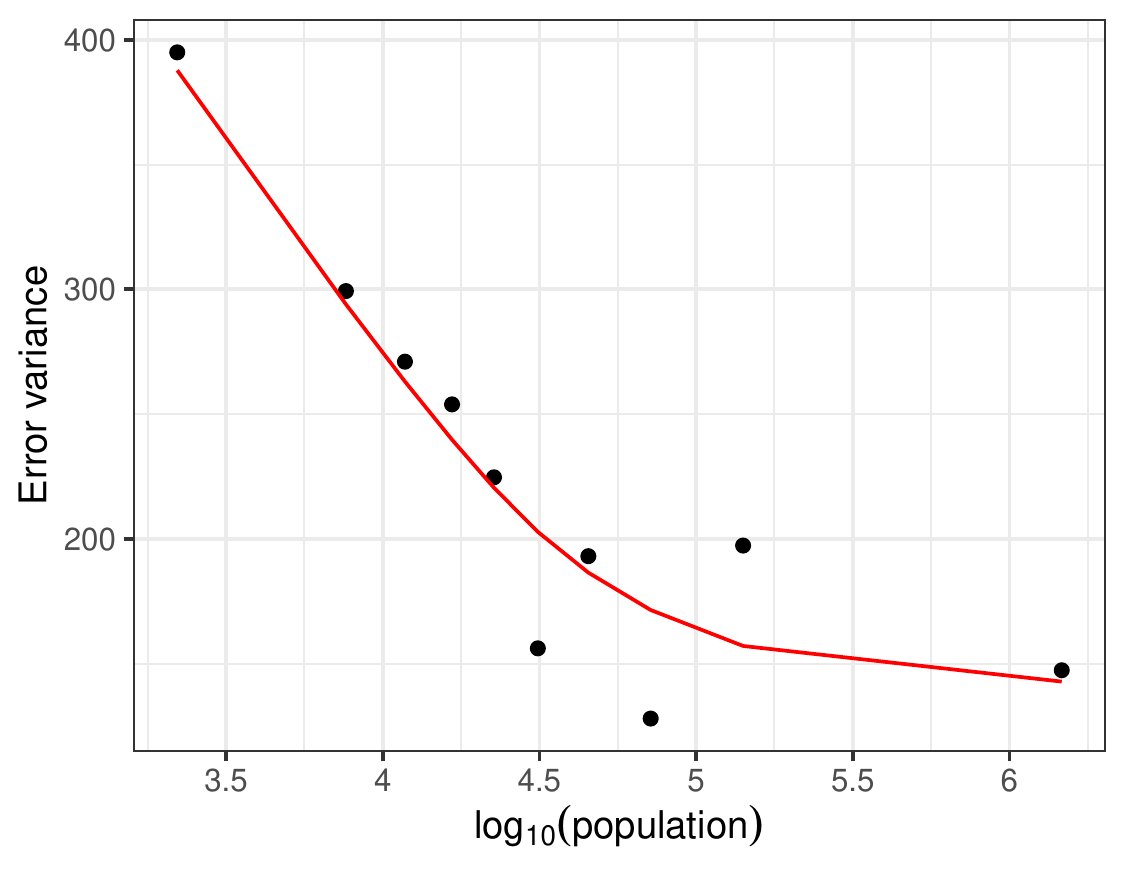}
	\caption{County-level cancer death rate of the states around the Great Lakes area (left) and error variance as a function of $\log_{10}(\text{population})$ (right).}
	\label{fig: rda_cancer_us_map}
	\vspace{2em}
	
	\includegraphics[scale = 0.6]{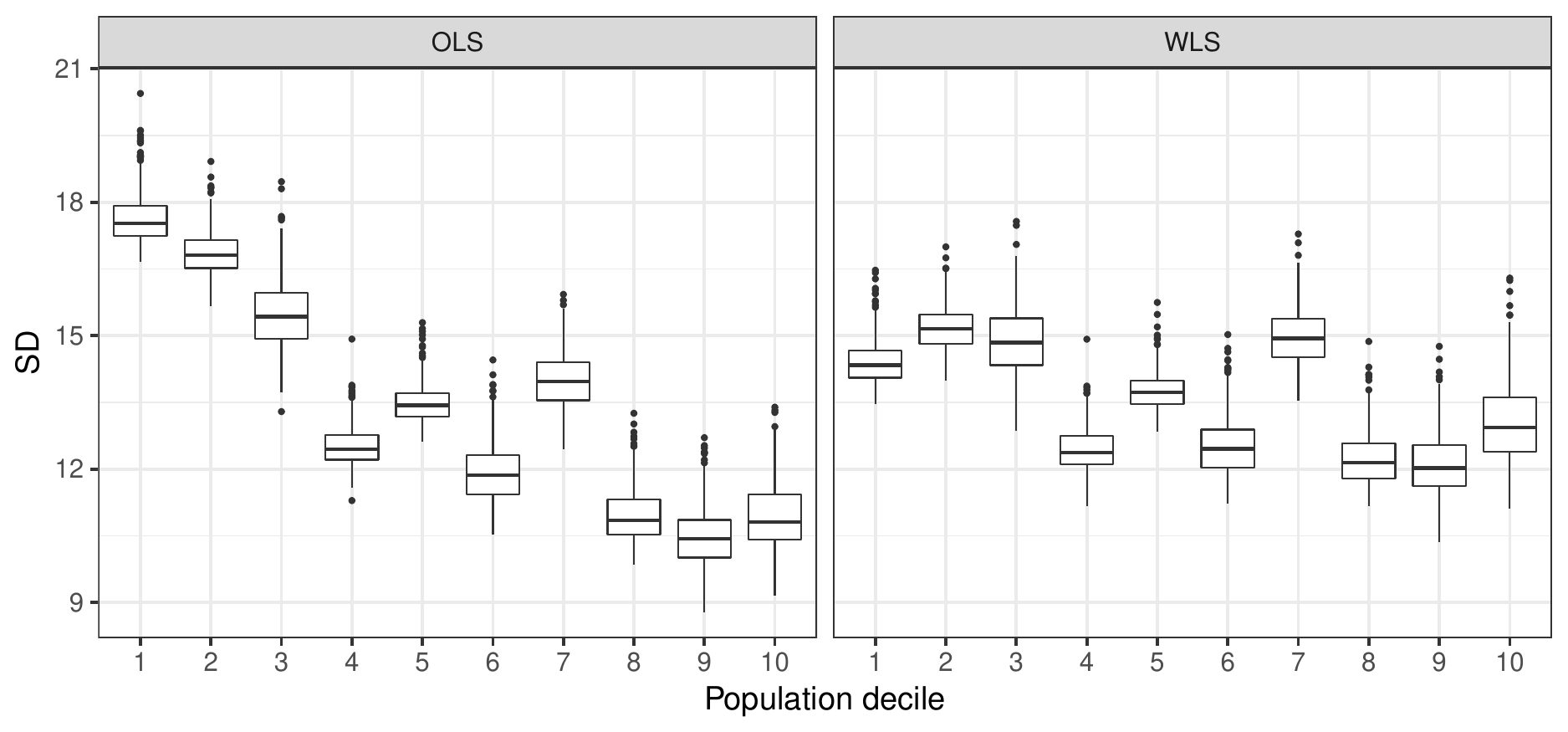}
    \caption{Standard deviation of the weighted prediction errors $\sqrt{w_i} (y_i - \hat{y}_i)$ by population deciles. For the OLS model, $w_i = 1$, whereas for the WLS model, $w_i$'s are calculated using the smoothing spline model shown in the right panel of Figure \ref{fig: rda_cancer_us_map}.}\label{fig: rda_cancer_weighted_errs_sd}
\end{figure}

We revisit the county-level cancer mortality data discussed in \cite{luan2021predictive}. In addition to the 22 continuous features considered in the previous analysis, we also take geographic variations into account by grouping all states into 4 regions, namely Midwest, Northeast, South and West, defined in Chapter 6 of the Geographic Areas Reference Manual (GARM) \citep{bureau1994geographic}. To restrict the scope of the analysis, we focus on the 19 states around the Great Lakes area illustrated in the left panel of Figure \ref{fig: rda_cancer_us_map}, which amount to 1,148 observations (counties) in total. Minnesota is removed due to data quality issues. Stratified random sampling based on each region is used to create two datasets of size 150 and 998 for training and testing respectively.
 
We use observations from the remaining states to estimate the variance function $\tau$ and determine the variable orderings in subset regression. For estimation of $\tau$, we first fit an OLS model to the death rate with all features included. We then split the residuals of the model into 10 quantile groups by the the population size and calculate the variances of the residuals within each group. We then estimate the variance function by fitting a smoothing spline model to these sample variances as a function of the group median of the logarithm of the population size, as illustrated in the right panel of Figure \ref{fig: rda_cancer_us_map}. For ordering of variables, we bootstrap data for 500 times. With each bootstrap sample, we construct a variable sequence by adding the variable that reduces the weighted residual sum of squares most until all variables are included. The variables are then ordered by their average ranks over all 500 bootstrap samples.
\begin{figure}[t!]
    \centering
    \includegraphics[scale = 0.55]{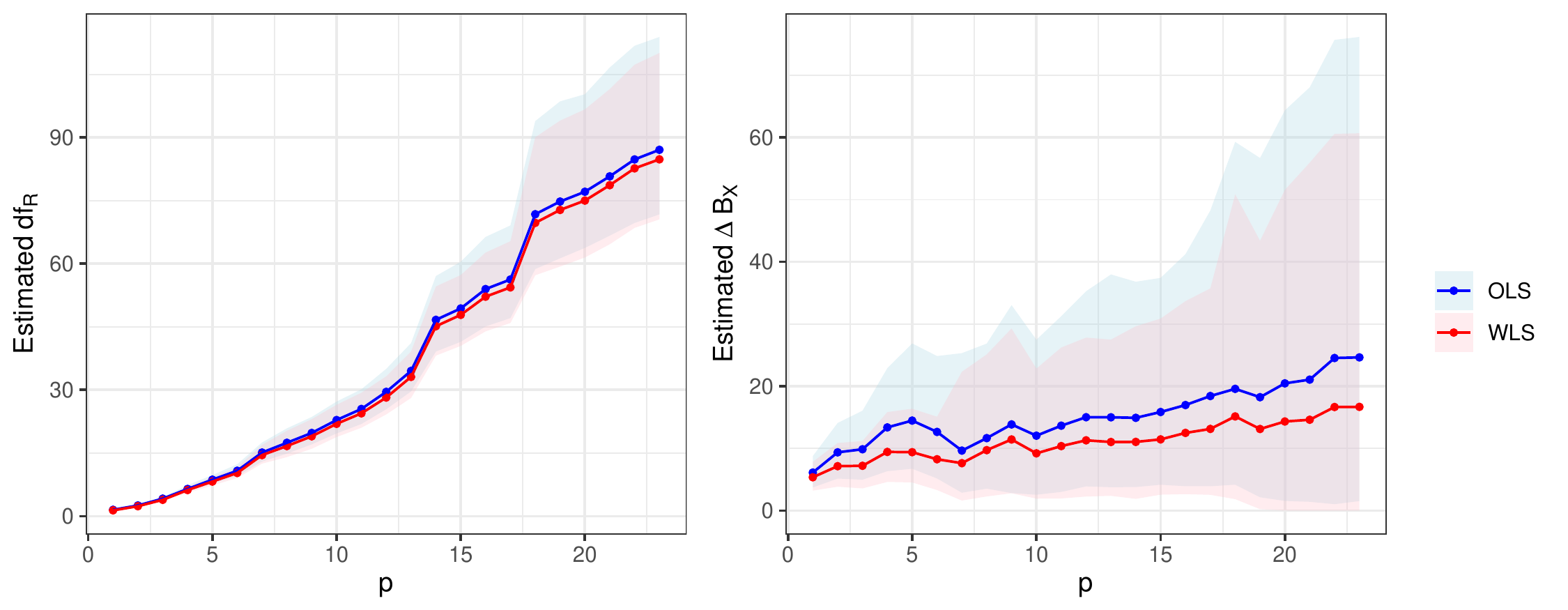}
    \caption{Estimated predictive model degrees of freedom $\mathrm{df}_{\rm R}$ and the excess bias $\Delta B_\bX$ of the ordinary and weighted least squares models averaged over 500 random splits of the data into training and test sets. The shaded bands show the 2.5th and 97.5th percentiles of the estimates for each $p$.}
    \label{fig: rda_cancer_df_comparison}
\end{figure}

We use $w = 1/\tau$ as the weight function in the weighted least squares model. We first compare the OLS and WLS models based on 500 random splits of the data into training and test sets. Figure \ref{fig: rda_cancer_weighted_errs_sd} shows the standard deviations of the weighted prediction errors $\sqrt{w_i} (y_i - \hat{y}_i)$ by population deciles. For the OLS model, these errors are the unweighted ones due to $w_i = 1$. It is clear that error variances become much more homogeneous when case weights are incorporated. In Figure \ref{fig: rda_cancer_df_comparison}, we demonstrate the average estimated predictive model degrees of freedom and excess bias of the two models. We can see that the WLS model has both smaller model complexity and excess bias than the OLS model, though the difference in the former is rather negligible. In terms of model selection, the two models are fairly close, with the OLS model choosing 9.07 (SD = 1.65) variables on average and 9.08 (SD = 1.78) for the WLS model. For a particular random split of the data, a comparison between the WLS and OLS models in coefficient estimation is given in Appendix \ref{appendix: rda_cancer}.

\begin{table}[!t]
	\centering
	\caption{Average performance of different variable selection methods based on 500 random data partitions. $p_\ast$ and $\hat{p}_\ast$ are the optimal model size indicated by the test data and variable selection criterion respectively. Numbers in parentheses are the standard deviations.}
	\label{tab: rda_cancer_subset_reg_avg}
	\begin{tabular}{c|cccccc}
		\multicolumn{7}{c}{$n = 150$}\\
		\hline
		Method & $\widehat{\wErrR}$ & $C_p$ & LOOCV & 5-fold CV & AIC & BIC \\ 
		\hline
		$\hat{p}_\ast/p_\ast$ & \makecell{$0.968$\\\small($0.450$)}& \makecell{$1.191$\\\small($0.453$)} & \makecell{$1.033$\\\small($0.458$)} & \makecell{$1.062$\\\small($0.460$)} &  \makecell{$1.247$\\\small($0.469$)} & \makecell{$0.703$\\\small($0.317$)}  \\
		
		$\wErrR_{\rm te}(\hat{p}_\ast)/\wErrR_{\rm te}(p_\ast)$ & \makecell{$1.071$\\\small($0.123$)}& \makecell{$1.071$\\\small($0.117$)} & \makecell{$1.107$\\\small($0.123$)} & \makecell{$1.064$\\\small($0.109$)} & \makecell{$1.073$\\\small($0.123$)} & \makecell{$1.085$\\\small($0.125$)}  \\
		\hline
		\multicolumn{7}{c}{$n = 40$}\\
		\hline
		Method & $\widehat{\wErrR}$ & $C_p$ & LOOCV & 5-fold CV & AIC & BIC \\ 
		\hline
		$\hat{p}_\ast/p_\ast$ & \makecell{$0.663$\\\small($0.319$)}& \makecell{$1.492$\\\small($0.707$)} & \makecell{$1.230$\\\small($0.616$)} & \makecell{$1.155$\\\small($0.554$)} &  \makecell{$2.146$\\\small($0.840$)} & \makecell{$1.192$\\\small($0.612$)}  \\
		
		$\wErrR_{\rm te}(\hat{p}_\ast)/\wErrR_{\rm te}(p_\ast)$ & \makecell{$1.370$\\\small($0.530$)}& \makecell{$1.936$\\\small($1.719$)} & \makecell{$1.689$\\\small($1.440$)} & \makecell{$1.604$\\\small($1.105$)} & \makecell{$2.589$\\\small($2.529$)} & \makecell{$1.703$\\\small($1.386$)}  \\
		\hline
	\end{tabular}
\end{table}
In Table \ref{tab: rda_cancer_subset_reg_avg}, we compare the average performance of $\widehat{\wErrR}$ in variable selection with five other commonly used criteria using the estimated weight scheme. Here, we consider two scenarios. In the first one, we set the sample size of the training data equal to 150. On average, models selected by $\widehat{\wErrR}$ and BIC contain fewer variables than the optimal ones on the test data. However, since the sample size is much larger than the total number of variables in this case, the model sizes identified by each method are not too different. In terms of out-of-sample prediction, all criteria have similar performance.

If we reduce the sample size of training data to 40, more differences can be observed among these six criteria. In general, $\widehat{\wErrR}$ outperforms others since it selects the least number of variables while having the smallest prediction error and lowest variability. The model parsimony from $\widehat{\wErrR}$ is due to the fact that model complexity increases quickly when the number of features gets close to the sample size. Since out-of-sample prediction in this scenario usually involves much higher uncertainty, it is reasonable to only include the most relevant variables in the model.

\subsubsection{King County House Sales Data} \label{subsubsec: rda_house}

We continue the discussion in Section \ref{subsubsec: wls_vs_ols} on the impact of weight scheme on model selection using the King County house sales data, which is publicly available at \url{https://www.kaggle.com/harlfoxem/housesalesprediction}. It contains sale prices and property features of 21,597 houses sold in King County, WA between May 2014 and May 2015. We split the data by year (2014 and 2015) for exploration and modeling purposes. We take square root of the housing price as the response variable. After selecting 12 initial predictors, we apply appropriate transformations. Further, using the same data, we estimate the variance function $\tau$ based on the residuals from the OLS model with all features included, as shown by the red curve in Figure \ref{fig: rda_house_variance_function}. Details of the procedure are provided in Appendix \ref{appendix: rda_house_sales}.
\begin{figure}[!t]
    \centering
    \includegraphics[scale = 0.55]{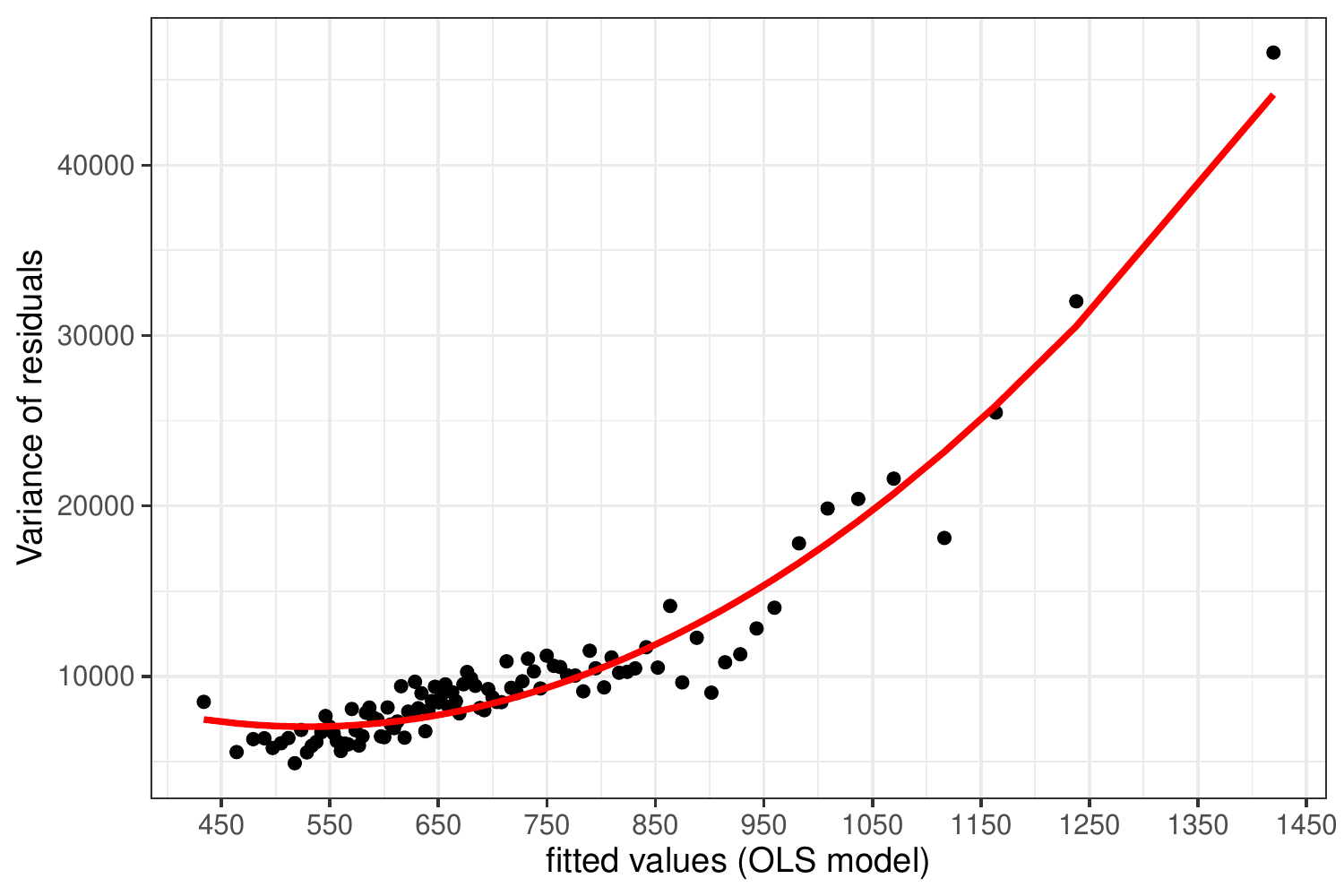}
    \caption{Variances of residuals within each percentile group of the fitted values of the OLS model with the square root of the housing price as a response and all 12 predictors for King County house sales data.} \label{fig: rda_house_variance_function}
\end{figure}

For the modeling part, we randomly select a subset of size $n = 200$ to fit a set of WLS models using the estimated variance function $\tau$ along the variable sequence. We repeat the procedure 500 times. For comparison, we also fit an OLS model with the same set of variables for each subset size $p$. We find that, under both in-sample and out-of-sample prediction settings, the two models identify the same subset size for more than half of the time, which is apparently different from what we have seen from the simulation study in Section \ref{subsubsec: wls_vs_ols}. This suggests that the training error and excess bias in this example may have larger influence on model selection. What is consistent with the simulation is that, when the difference does show up, it is most likely that the WLS model selects more variables than the OLS model. In a few extreme cases, this difference in the size of the selected model can reach up to 6 variables. 

For a close-up examination of the above findings, we decompose $\widehat{\wErrF}$ and $\widehat{\wErrR}$ into the training error, excess bias and excess variance, and study the median percentages of these components of the selected models. A figure of the decomposition is provided in Appendix \ref{appendix: rda_house_sales}. We find that the training error is the dominant component in the risk estimate under both in-sample and out-of-sample prediction settings, even when all variables are included. As a result, model selection is greatly impacted by the training error, especially in the in-sample scenario where the term explains more than 90\% of the risk. Since the training error is decreasing in the subset size, both WLS and OLS models tend to favor the full model in this scenario. On the other hand, for out-of-sample prediction, the excess variance plays a relatively more important role in model selection (around 20-25\% of the risk when all variables are included against 9-12\% in the in-sample prediction setting), which leads to a reduction in the size of the selected model. Comparing the WLS and OLS models, the former has a constantly higher percentage of training error and lower percentage of excess variance than the latter. This explains why the WLS model has a slight tendency of selecting a few more variables. In addition, we also see a lower percentage of the excess bias in $\widehat{\wErrR}$ of the WLS model. Since the overall percentage of this component is quite small in this example, we expect its impact on model selection would be small as well. 

\section{Conclusion} \label{sec: conclusion}

In this work, we extend the existing theories on model complexity to a heteroscedastic setting for linear regression procedures. In particular, we generalize the classical model degrees of freedom and predictive model degrees of freedom so that models that incorporate different case weights can be distinguished appropriately. As a typical example, we explore the two extended measures of model degrees of freedom of the weighted least squares method in depth. Interestingly, we find that, in the subset regression setting, the two measures are a shrunken version of those of the ordinary least squares method, and the amount of shrinkage depends on how uneven the error variances are. Furthermore, we set up connections between the two measures of model degrees of freedom and weight schemes both for model fitting and model evaluation. Our results provide valuable insights on how model complexity is affected by case weights, and further provide useful guidance in risk estimation and model selection.

For future study in model complexity, it will be interesting to extend the current framework to a generalized linear model setting. Note that for exponential family distributions other than normal, variance usually shares the same set of parameters with the mean and varies from sample to sample. As a result, our work on heteroscedastic linear regression models holds promise for the possibility of the extension and thus deserves further investigation.

\subsubsection*{Acknowledgements}

We thank Steven MacEachern for comments on the earlier version of this work. This research was supported in part by the National Science Foundation Grants DMS-17-12580, DMS-17-21445 and DMS-20-15490.

\bibliographystyle{apalike}
\bibliography{bibliography}

\begin{thebibliography}{}

\bibitem[Anscombe, 1961]{anscombe1961examination}
Anscombe, F.~J. (1961).
\newblock Examination of residuals.
\newblock In {\em Proceedings of the Fourth Berkeley Symposium on Mathematical
  Statistics and Probability, Volume 1: Contributions to the Theory of
  Statistics}, pages 1--36. University of California Press.

\bibitem[Bickel, 1978]{bickel1978using}
Bickel, P.~J. (1978).
\newblock Using residuals robustly {I}: Tests for heteroscedasticity,
  nonlinearity.
\newblock {\em Annals of Statistics}, 6(2):266--291.

\bibitem[Box and Hill, 1974]{box1974correcting}
Box, G.~E. and Hill, W.~J. (1974).
\newblock Correcting inhomogeneity of variance with power transformation
  weighting.
\newblock {\em Technometrics}, 16(3):385--389.

\bibitem[Breiman et~al., 2017]{breiman2017classification}
Breiman, L., Friedman, J.~H., Olshen, R.~A., and Stone, C.~J. (2017).
\newblock {\em {Classification and Regression Trees}}.
\newblock Routledge.

\bibitem[{Bureau of the Census}, 1994]{bureau1994geographic}
{Bureau of the Census} (1994).
\newblock {\em Geographic Areas Reference Manual}.
\newblock U.S. Department of Commerce, Economics and Statistics Administration.

\bibitem[Carroll and Ruppert, 1981]{carroll1981robust}
Carroll, R.~J. and Ruppert, D. (1981).
\newblock On robust tests for heteroscedasticity.
\newblock {\em Annals of Statistics}, 9(1):206--210.

\bibitem[Craven and Wahba, 1978]{craven1978smoothing}
Craven, P. and Wahba, G. (1978).
\newblock Smoothing noisy data with spline functions.
\newblock {\em Numerische {Mathematik}}, 31(4):377--403.

\bibitem[Dempster et~al., 1977]{dempster1977maximum}
Dempster, A.~P., Laird, N.~M., and Rubin, D.~B. (1977).
\newblock Maximum likelihood from incomplete data via the {EM} algorithm.
\newblock {\em Journal of the Royal Statistical Society: Series B
  (Methodological)}, 39(1):1--22.

\bibitem[Efron, 1986]{efron1986biased}
Efron, B. (1986).
\newblock How biased is the apparent error rate of a prediction rule?
\newblock {\em Journal of the American Statistical Association},
  81(394):461--470.

\bibitem[Efron, 2004]{efron2004estimation}
Efron, B. (2004).
\newblock The estimation of prediction error: Covariance penalties and
  cross-validation.
\newblock {\em Journal of the American Statistical Association},
  99(467):619--632.

\bibitem[Lowd and Domingos, 2005]{lowd2005naive}
Lowd, D. and Domingos, P. (2005).
\newblock Naive {Bayes} models for probability estimation.
\newblock In {\em Proceedings of the 22nd International Conference on Machine
  Learning}, pages 529--536.

\bibitem[Luan et~al., 2021]{luan2021predictive}
Luan, B., Lee, Y., and Zhu, Y. (2021).
\newblock Predictive model degrees of freedom in linear regression.
\newblock {\em arXiv preprint arXiv:2106.15682}.

\bibitem[Mallows, 1973]{mallows1973some}
Mallows, C. (1973).
\newblock Some comments on {$C_p$}.
\newblock {\em Technometrics}, 15(4):661--675.

\bibitem[Mardia et~al., 1979]{mardia1979multivariate}
Mardia, K., Kent, J., and Bibby, J. (1979).
\newblock {\em {Multivariate Analysis}}.
\newblock Academic Press.

\bibitem[Pivaro et~al., 2017]{pivaro2017exact}
Pivaro, G.~F., Kumar, S., Fraidenraich, G., and Dias, C.~F. (2017).
\newblock On the exact and approximate eigenvalue distribution for sum of
  {Wishart} matrices.
\newblock {\em IEEE Transactions on Vehicular Technology}, 66(11):10537--10541.

\bibitem[Reiter, 2005]{reiter2005using}
Reiter, J.~P. (2005).
\newblock Using {CART} to generate partially synthetic public use microdata.
\newblock {\em Journal of Official Statistics}, 21(3):441.

\bibitem[Tibshirani and Hastie, 1987]{tibshirani1987local}
Tibshirani, R. and Hastie, T. (1987).
\newblock Local likelihood estimation.
\newblock {\em Journal of the American Statistical Association},
  82(398):559--567.

\bibitem[Ye, 1998]{ye1998measuring}
Ye, J. (1998).
\newblock On measuring and correcting the effects of data mining and model
  selection.
\newblock {\em Journal of the American Statistical Association},
  93(441):120--131.

\end{thebibliography}

\newpage

\begin{appendices}

\section{Proofs}

\subsection{Derivation of Equations \eqref{eq: bias_var_ErrR_weighted} and \eqref{eq: bias_var_ErrT_weighted}} \label{pf: bias_var_Err_weighted}

Consider a linear procedure with hat matrix $\bH$ and hat vector $\bh_\ast$ at $\bx_\ast$. We first look at the bias-variance decomposition for $\wErrR_\bX$. Note that
$$
    \wErrR_\bX = \E[w_\ast (y_\ast - \bh_\ast^\T \by)^2 \vert \bX] = \E[\E[(w_\ast (y_\ast - \bh_\ast^\T \by)^2 \vert \bx_\ast, \bX]\vert \bX],
$$
where the inner and outer expectations are taken with respect to $(y_\ast, \by)$ and $\bx_\ast$ respectively. For the inner expectation, using the fact that $(\bx_\ast,y_\ast)$ is independent of $(\bX,\by)$, we have
$$
\begin{aligned}
    \E[(w_\ast (y_\ast - \bh_\ast^\T \by)^2 \vert \bx_\ast, \bX]
    & = w_\ast \E[(y_\ast - \mu_\ast + \mu_\ast - \bh_\ast^\T\bmu + \bh_\ast^\T\bmu - \bh_\ast^\T \by)^2 \vert \bx_\ast, \bX] \\
    & = w_\ast(\E[(y_\ast-\mu_\ast)^2\vert\bx_\ast] + (\mu_\ast - \bh_\ast^\T\bmu)^2 + \E[(\bh_\ast^\T\by - \bh_\ast^\T \bmu)^2\vert \bx_\ast,\bX])\\
    & = \sigma_\varepsilon^2 w_\ast \tau_\ast + w_\ast (\mu_\ast - \bh_\ast^\T\bmu)^2 + \sigma_\varepsilon^2 w_\ast \Vert \bh_\ast \Vert_\bT^2,
\end{aligned}
$$
where $\bT = \mathrm{diag}(\tau_1,\ldots,\tau_n)$. Then taking another expectation of the above result with respect to $\bx_\ast$ yields \eqref{eq: bias_var_ErrR_weighted} immediately.

For $\ErrT_\bX$, it follows that
$$
\begin{aligned}
    \wErrT_\bX 
    & = \frac{1}{n}\E(\Vert \by - \bH \by \Vert_\bW^2 \vert \bX)\\
    & = \frac{1}{n} \Vert \bmu-\bH\bmu \Vert_\bW^2 + \frac{1}{n}\sigma_\varepsilon^2 \trace[(\bI_n-\bH)^\T \bW (\bI_n-\bH)\bT] \\
    & = \frac{1}{n}\sigma_\varepsilon^2 \trace(\bW\bT) + \frac{1}{n} \Vert \bmu-\bH\bmu \Vert_\bW^2 + \frac{1}{n}\sigma_\varepsilon^2 \trace(\bH^\T \bW \bH \bT - \bH^\T \bW \bT - \bW \bH \bT)\\
    & = \frac{1}{n}\sigma_\varepsilon^2 \sum_{i=1}^{n}{w_i\tau_i} + \frac{1}{n} \Vert \bmu-\bH\bmu \Vert_\bW^2 + \frac{1}{n}\sigma_\varepsilon^2 \trace(\bH^\T \bW \bH \bT - 2\bH \bT \bW).
\end{aligned}
$$
This completes the proof of \eqref{eq: bias_var_ErrT_weighted}. In particular, when $w(\cdot) = \tau(\cdot) = 1$, \eqref{eq: bias_var_ErrR_weighted} and \eqref{eq: bias_var_ErrT_weighted} reduce to \eqref{eq: bias_var_ErrR_unweighted} and \eqref{eq: bias_var_ErrT_unweighted} respectively.

\subsection{Proof of Proposition \ref{prop: dfR_cov_representation}}

\begin{proof}
    Since $\hat{\mu}_\ast = \bh_\ast^\T \by$ and $\hat{\mu}_j = \bh_j^\T \by$, we have
    $$
    	\frac{\cov(y_i, \hat{\mu}_\ast \vert \bx_\ast, \bX)}{\sigma_\varepsilon^2 \tau_i} = \frac{\partial \E(\hat{\mu}_\ast\vert \bx_\ast,\bX)}{\partial \mu_i} = h_{\ast,i} \text{ and } \frac{\cov(y_i, \hat{\mu}_j \vert \bx_\ast, \bX)}{\sigma_\varepsilon^2 \tau_i} = \frac{\partial \E(\hat{\mu}_j\vert \bX)}{\partial \mu_i} = h_{ji}.
    $$
    As a result,
    \begin{small}
    $$
    	\sum_{i=1}^n \tau_i \E\left(w_\ast \frac{\cov^2(y_i, \hat{\mu}_\ast \vert \bx_\ast, \bX)}{(\sigma_\varepsilon^2 \tau_i)^2} \Bigg\vert \bX\right) = \sum_{i=1}^n \tau_i \E\left(w_\ast \left(\frac{\partial \E(\hat{\mu}_\ast \vert \bx_\ast, \bX) }{\partial \mu_i}\right)^2 \Bigg\vert \bX\right) = \E(w_\ast \Vert \bh_\ast \Vert_\bT^2 \vert \bX),
    $$
    \end{small}
    and
    \begin{small}
    $$
    	\sum_{i=1}^n \tau_i \sum_{j=1}^n w_j \frac{\cov^2(y_i, \hat{\mu}_j \vert \bX)}{(\sigma_\varepsilon^2 \tau_i)^2} = \sum_{i=1}^n \tau_i \sum_{j=1}^n w_j \left(\frac{\partial \E(\hat{\mu}_j \vert \bX) }{\partial \mu_i}\right)^2 = \trace(\bW \bH \bT \bH^\T).
    $$
    \end{small}
    The theorem then follows.
\end{proof}

\subsection{Proof of Theorem \ref{thm: DdfF_wls_general}} \label{pf: DdfF_wls_general}

\begin{proof}
    Let $\bU =  \frac{1}{\bar{w}} \bW^{\frac{1}{2}}\bT \bW^{\frac{1}{2}}$. We can express $\dfF(p)$ in terms of $\tilde{\bX}_p$ and $\bU$ as
    $$
        \dfF(p) = \trace[(\tilde{\bX}_p^\T \tilde{\bX}_p)^{-1} \tilde{\bX}_p^\T \bU \tilde{\bX}_p].
    $$
    Note that $\tilde{\bX}_{p+1} = (\tilde{\bX}_p, \tilde{\bx}_{(p+1)})$. By the inversion formula for a 2 by 2 partitioned matrix,
    $$
        (\tilde{\bX}_{p+1}^\T \tilde{\bX}_{p+1})^{-1} = 
        \begin{pmatrix}
            \tilde{\bX}_p^\T \tilde{\bX}_p & \tilde{\bX}_p^\T \tilde{\bx}_{(p+1)}\\
            \tilde{\bx}_{(p+1)}^\T \tilde{\bX}_p & \tilde{\bx}_{(p+1)}^\T \tilde{\bx}_{(p+1)}
        \end{pmatrix}^{-1}
        = \left(
        \begin{array}{cc}
        	(\tilde{\bX}_p^\T \tilde{\bX}_p)^{-1} + \frac{1}{b} \bs \bs^\T & -\frac{1}{b} \bs \\
        	-\frac{1}{b} \bs^\T  & \frac{1}{b}
        \end{array}\right),
    $$
    where $\bs = (\tilde{\bX}_p^\T \tilde{\bX}_p)^{-1} \tilde{\bX}_p^\T \tilde{\bx}_{(p+1)}$ and $b = \tilde{\bx}_{(p+1)}^\T (\bI_n - \tilde{\bX}_p (\tilde{\bX}_p^\T \tilde{\bX}_p)^{-1} \tilde{\bX}_p^\T) \tilde{\bx}_{(p+1)}$. Letting $\tilde{\bH} = \tilde{\bX}_p (\tilde{\bX}_p^\T \tilde{\bX}_p)^{-1} \tilde{\bX}_p^\T$ and $\tilde{\br}_{(p+1)} = (\bI_n - \tilde{\bH})\tilde{\bx}_{(p+1)}$, we have $b = \Vert \tilde{\br}_{(p+1)} \Vert^2$. Then
    $$
    \begin{aligned}
        \dfF(p+1) 
        & = \trace[(\tilde{\bX}_{p+1}^\T \tilde{\bX}_{p+1})^{-1} \tilde{\bX}_{p+1}^\T \bU \tilde{\bX}_{p+1}] \\
        & = \trace\left[\left((\tilde{\bX}_p^\T \tilde{\bX}_p)^{-1} + \frac{1}{b}\bs\bs^\T \right) \tilde{\bX}_p^\T \bU \tilde{\bX}_p - \frac{1}{b} \bs \tilde{\bx}_{(p+1)}^\T \bU \tilde{\bX}_p \right] \\
        & \qquad - \frac{1}{b} \bs^\T \tilde{\bX}_p^\T \bU \tilde{\bx}_{(p+1)} + \frac{1}{b} \tilde{\bx}_{(p+1)}^\T \bU \tilde{\bx}_{(p+1)}\\
        & = \dfF(p) + \frac{1}{b} \tilde{\bx}_{(p+1)}^\T \tilde{\bX}_p (\tilde{\bX}_p^\T \tilde{\bX}_p)^{-1} \tilde{\bX}_p^\T \bU \tilde{\bX}_p (\tilde{\bX}_p^\T \tilde{\bX}_p)^{-1} \tilde{\bX}_p^\T \tilde{\bx}_{(p+1)}\\
        & \qquad\qquad\, - \frac{2}{b} \tilde{\bx}_{(p+1)}^\T \tilde{\bX}_p (\tilde{\bX}_p^\T \tilde{\bX}_p)^{-1} \tilde{\bX}_p^\T \bU \tilde{\bx}_{(p+1)}\\
        & \qquad\qquad\, + \frac{1}{b} \tilde{\bx}_{(p+1)}^\T \bU \tilde{\bx}_{(p+1)}\\
        & = \dfF(p) + \frac{1}{b} \tilde{\bx}_{(p+1)}^\T (\bI_n - \tilde{\bH}) \bU (\bI_n - \tilde{\bH}) \tilde{\bx}_{(p+1)} \\
        & = \dfF(p) + \frac{1}{b} \tilde{\br}_{(p+1)}^\T \bU \tilde{\br}_{(p+1)}\\
        & = \dfF(p) + \frac{\Vert \tilde{\br}_{(p+1)} \Vert_\bU^2}{\Vert \tilde{\br}_{(p+1)} \Vert^2}.
    \end{aligned}
    $$
    
\end{proof}

\subsection{Proof of Theorem \ref{thm: DdfR_wls_general}}
\begin{proof}
    Following the notation in \ref{pf: DdfF_wls_general}, define $\bG_p = (\tilde{\bX}_p^\T \tilde{\bX}_p)^{-1} \tilde{\bX}_p^\T \bU \tilde{\bX}_p (\tilde{\bX}_p^\T \tilde{\bX}_p)^{-1}$. Then $\dfR$ can be written as
    $$
        \dfR(p) = \frac{1}{2}\dfF(p) + \frac{n}{2} \trace(\bG_p \bSigma_p).
    $$
    It can be checked that $\bG_{p+1}=(\tilde{\bX}_{p+1}^\T \tilde{\bX}_{p+1})^{-1} \tilde{\bX}_{p+1}^\T \bU \tilde{\bX}_{p+1} (\tilde{\bX}_{p+1}^\T \tilde{\bX}_{p+1})^{-1}$ has the form
    $$
        \bG_{p+1} =
        \begin{pmatrix}
            \bF & \bg\\
            \bg^\T & f
        \end{pmatrix},
    $$
    where
    $$
    \begin{aligned}
        \bF & = \bG_p - \frac{1}{b} (\tilde{\bX}_p^\T \tilde{\bX}_p)^{-1} \tilde{\bX}_p^\T \bU \tilde{\br}_{(p+1)} \bs^\T - \frac{1}{b} \bs \tilde{\br}_{(p+1)}^\T \bU \tilde{\bX}_p (\tilde{\bX}_p^\T \tilde{\bX}_p)^{-1} + \frac{\Vert \tilde{\br}_{(p+1)} \Vert_{\bU}^{2}}{b^2} \bs \bs^\T, \\
        \bg & =  \frac{1}{b} (\tilde{\bX}_p^\T \tilde{\bX}_p)^{-1} \tilde{\bX}_p^\T \bU \tilde{\br}_{(p+1)} - \frac{\Vert \tilde{\br}_{(p+1)} \Vert_{\bU}^{2}}{b^2} \bs,\\
        f & = \frac{\Vert \tilde{\br}_{(p+1)} \Vert_{\bU}^{2}}{b^2}.
    \end{aligned}
    $$
    Note that
    $$
        \bF = \bG_p-(\bg+f\bs)\bs^\T-\bs(\bg+f\bs)^\T+f\bs\bs^\T = \bG_p - \bg\bs^\T - \bs\bg^\T - f\bs\bs^\T.
    $$ 
    Then
    $$
    \begin{aligned}
        \dfR(p+1)
        & = \frac{1}{2}\dfF(p+1) + \frac{n}{2} \trace(\bG_{p+1} \bSigma_{p+1}) \\
        & = \frac{1}{2}\dfF(p+1) + \frac{n}{2} [\trace(\bF \bSigma_p) + 2 \bg^\T \bphi_{p,p+1} + f \sigma_{p+1}^2 ]\\
        & = \dfR(p) + \frac{1}{2}\Delta\dfF(p) + \frac{n}{2}(f\sigma_{p+1}^2 - f\bs^\T \bSigma_p \bs - 2\bg^\T \bSigma_p \bs + 2\bg^\T \bphi_{p,p+1}).\\
    \end{aligned}
    $$
    Now we show that
    \begin{equation} \label{eq: DdfR_wls_general_equivalent_form}
    \begin{aligned}
        & \quad\; \frac{n}{2}(f\sigma_{p+1}^2 - f\bs^\T \bSigma_p \bs - 2\bg^\T \bSigma_p \bs + 2\bg^\T \bphi_{p,p+1}) \\
        & = \frac{n}{2} \Delta\dfF(p) \frac{(\Vert \tilde{\bx}_{(p+1)} - \tilde{\bX}_p \bSigma_p^{-1} \bphi_{p,p+1} \Vert_\mathbf{M}^2 + a_{p,p+1})}{\Vert \tilde{\br}_{(p+1)} \Vert^2} \\
        & \qquad\qquad - n \frac{\langle \tilde{\br}_{(p+1)}, \tilde{\bx}_{(p+1)} - \tilde{\bX}_p \bSigma_p^{-1} \bphi_{p,p+1} \rangle_{\bU\mathbf{M}}}{\Vert \tilde{\br}_{(p+1)} \Vert^2}.
    \end{aligned}
    \end{equation}
    
    Consider expressing the right hand side of \eqref{eq: DdfR_wls_general_equivalent_form} using $f$, $\bg$ and $\bs$. For the first term, note that
    $$
    \begin{aligned}
        \Vert \tilde{\bx}_{(p+1)} \Vert_\bM^2 & = \bs^\T \bSigma_p \bs, \\
        \Vert \tilde{\bX}_p \bSigma_p^{-1} \bphi_{p,p+1} \Vert_\bM^2 & = \bphi_{p,p+1}^\T \bSigma_p^{-1} \bphi_{p,p+1},\\
        \langle \tilde{\bx}_{(p+1)}, \tilde{\bX}_p \bSigma_p^{-1} \bphi_{p,p+1}\rangle_\bM & = \bs^\T \bphi_{p,p+1}.
    \end{aligned}
    $$
    Then 
    \begin{equation} \label{eq: DdfR_wls_general_equivalent_form_RHS1}
    \begin{aligned}
        & \quad\; \frac{n}{2} \Delta\dfF(p) \frac{(\Vert \tilde{\bx}_{(p+1)} - \tilde{\bX}_p \bSigma_p^{-1} \bphi_{p,p+1} \Vert_\mathbf{M}^2 + a_{p,p+1})}{\Vert \tilde{\br}_{(p+1)} \Vert^2}\\
        & = \frac{nf}{2} (\Vert \tilde{\bx}_{(p+1)} \Vert_\bM^2 - 2 \langle \tilde{\bx}_{(p+1)}, \tilde{\bX}_p \bSigma_p^{-1} \bphi_{p,p+1}\rangle_\bM + \Vert \tilde{\bX}_p \bSigma_p^{-1} \bphi_{p,p+1} \Vert_\bM^2 + a_{p,p+1})\\
        & = \frac{nf}{2}(\bs^\T \bSigma_p \bs - 2\bs^\T \bphi_{p,p+1} + \sigma_{p+1}^2).
    \end{aligned}
    \end{equation}
    For the second term, it can be similarly shown that
    \begin{equation} \label{eq: DdfR_wls_general_equivalent_form_RHS2}
        n \frac{\langle \tilde{\br}_{(p+1)}, \tilde{\bx}_{(p+1)} - \tilde{\bX}_p \bSigma_p^{-1} \bphi_{p,p+1} \rangle_{\bU\mathbf{M}}}{\Vert \tilde{\br}_{(p+1)} \Vert^2} = n(\bg + f\bs)^\T (\bSigma_p \bs - \bphi_{p,p+1})
    \end{equation}
    Subtracting \eqref{eq: DdfR_wls_general_equivalent_form_RHS2} from \eqref{eq: DdfR_wls_general_equivalent_form_RHS1} yields \eqref{eq: DdfR_wls_general_equivalent_form} immediately.
\end{proof}

\subsection{Derivation of Approximation \eqref{eq: dfR_wls_normal_approx}} \label{pf: dfR_wls_normal_approx}
Assume that $\bx_\ast$ is multivariate normal with $\E(\bx_\ast) = \mathbf{0}$ and $\var(\bx_\ast) = \bSigma_\bx$. Let $\bu_1,\ldots,\bu_n$ be i.i.d.\ random variables that are independent of $\bx_1,\ldots,\bx_n$. Define $w_i = 1/\tau(\bu_i)$ and $w_\ast = 1/\tau(\bu_\ast)$. Then, we have
$$
    \bSigma_\bx = \frac{\bar{w}}{\E(w_\ast)} \frac{1}{\bar{w}}\E(w_\ast \bx_\ast \bx_\ast^\T) = \frac{\bar{w}}{\E(w_\ast)} \bSigma_p.
$$
Let $\bz_i = \bSigma_\bx^{-\frac{1}{2}}\bx_i$ and $\bZ = (\bz_1,\ldots,\bz_n)^\T$. Note that $\bz_i$'s are i.i.d.\ $\mathcal{N}(\mathbf{0},\bI_p)$. Thus, $\bz_i \bz_i^\T$'s are independent Wishart random matrices with degrees of freedom 1 and scale matrix $\bI_p$, denoted by $\mathcal{W}(1,\bI_p)$. Note that, when $w \propto 1/\tau$,
$$
    \dfR(p) = \frac{\zeta_\tau}{2} (p + n \bar{w} \trace[(\bX_p^\T \bW \bX_p)^{-1} \bSigma_p)])
$$
and
$$
    n\bar{w} \trace[(\bX^\T \bW \bX)^{-1} \bSigma_p)] = \frac{\E(w_\ast)}{\bar{w}} n\bar{w} \trace[(\bZ^\T \bW \bZ)^{-1}] = \frac{\E(w_\ast)}{\bar{w}} \trace\left[\left(\sum_{i=1}^n \frac{w_i}{n\bar{w}} \bz_i \bz_i^\T\right)^{-1}\right].
$$
When $n$ is large, \cite{pivaro2017exact} suggest to approximate $\sum_{i=1}^n \frac{w_i}{n \bar{w}} \bz_i \bz_i^\T$ by
$$
    \sum_{i=1}^n \frac{w_i}{n\bar{w}} \bz_i \bz_i^\T \stackrel{\rm d}{\approx} \frac{\bV}{n_\tau}, \quad n_\tau = \left[\frac{1}{\sum_{i=1}^n \left(\frac{w_i}{n\bar{w}}\right)^2}\right] = \left[\frac{\left(\sum_{i=1}^n 1/\tau_i\right)^2}{\sum_{i=1}^n (1/\tau_i)^2}\right],
$$
where $\bV \sim \mathcal{W}(n_\tau, \bI_p)$. Since $\frac{\E(w_\ast)}{\bar{w}} \approx 1$, we then have
$$
    n\bar{w} \left(\sum_{i=1}^n w_i \bz_i \bz_i^\T\right)^{-1} \stackrel{\rm d}{\approx} n_\tau \bV^{-1}.
$$
Assume $n_\tau > p+1$. Note that $\bV^{-1}$ follows an inverse-Wishart distribution with degrees of freedom $n_\tau$ and scale matrix $\bI_p$. It then follows that $\E[\trace(\bV^{-1})] = \frac{p}{n_\tau - p - 1}$ \citep{mardia1979multivariate}. Using $\frac{n_\tau p}{n_\tau - p - 1}$ as an approximation of $n\bar{w} \trace[(\bX^\T \bW \bX)^{-1} \bSigma_p)]$, we then have
$$
    \dfR(p) \approx \frac{p}{2} \zeta_\tau \left(1 + \frac{n_\tau}{n_\tau - p - 1}\right).
$$
\hfill\qedsymbol

\subsection{Proof of Theorem \ref{eq: differential_dfR_wls_general}} \label{pf: differential_dfR_wls_general}
\begin{proof}
    We first assume that $\bQ$ is an arbitrary $n \times n$ matrix such that $(\bX^\T \bQ \bX)^{-1}$ exists. Define $\bm{\Pi} \colon \R^{n \times n} \to \R^{p \times n}$ as
    $$
        \bm{\Pi}(\bQ) = (\bX^\T \bQ \bX)^{-1} \bX^\T \bQ.
    $$
    Then $\dfR$ can be written as
    $$
        \dfR = \zeta_\tau p + \frac{n}{2}\trace[\bm{\Pi}(\bQ) \bT \bm{\Pi}(\bQ)^\T \bR].
    $$
    Note that
    $$
    \begin{aligned}
        \partial \bm{\Pi}(\bQ) 
        & = -(\bX^\T \bQ \bX)^{-1} \bX^\T (\partial\bQ) \bX (\bX^\T \bQ \bX)^{-1} \bX^\T \bQ  + (\bX^\T \bQ \bX)^{-1} \bX^\T (\partial\bQ)\\
        & = (\bX^\T \bQ \bX)^{-1} \bX^\T (\partial\bQ) (\bI_n - \bH).
    \end{aligned}
    $$
    Then we have
    $$
    \begin{aligned}
        \partial(\bm{\Pi}(\bQ) \bT \bm{\Pi}(\bQ)^\T \bR) 
        & = (\partial \bm{\Pi}(\bQ)) \bT \bm{\Pi}(\bQ)^\T \bR + \bm{\Pi}(\bQ) \bT (\partial \bm{\Pi}(\bQ))^\T \bR\\
        & = (\bX^\T \bQ \bX)^{-1} \bX^\T (\partial\bQ) (\bI_n - \bH) \bT \bQ \bX (\bX^\T \bQ \bX)^{-1} \bR \\
        & \qquad + (\bX^\T \bQ \bX)^{-1} \bX^\T \bQ \bT (\bI_n - \bH)^\T (\partial\bQ) \bX (\bX^\T \bQ \bX)^{-1} \bR.
    \end{aligned}
    $$
    By the chain rule, it follows that
    $$
    \begin{aligned}
        \partial \dfR 
        & = \frac{n}{2} \partial \trace[\bm{\Pi}(\bQ) \bT \bm{\Pi}(\bQ)^\T \bR] \\
        & = \frac{n}{2} \trace[\partial(\bm{\Pi}(\bQ) \bT \bm{\Pi}(\bQ)^\T \bR)]\\
        & = \frac{n}{2} \trace[(\bX^\T \bQ \bX)^{-1} \bX^\T (\partial\bQ) (\bI_n - \bH) \bT \bQ \bX (\bX^\T \bQ \bX)^{-1} \bR ]\\
        & \qquad + \frac{n}{2} \trace[ (\bX^\T \bQ \bX)^{-1} \bX^\T \bQ \bT (\bI_n - \bH)^\T (\partial\bQ) \bX (\bX^\T \bQ \bX)^{-1} \bR ] \\
        & = n \trace[(\bI_n - \bH) \bT \bQ \bX (\bX^\T \bQ \bX)^{-1} \bR (\bX^\T \bQ \bX)^{-1} \bX^\T \partial\bQ ].
    \end{aligned}
    $$
    Thus,
    $$
        \frac{\partial \dfR}{\partial \bQ} = n (\bI_n - \bH) \bT \bQ \bX (\bX^\T \bQ \bX)^{-1} \bR (\bX^\T \bQ \bX)^{-1} \bX^\T.
    $$
    Since $\bQ$ is a diagonal matrix, we further have
    $$
        \frac{\partial \dfR}{\partial \bq} = n \, \mathrm{diag}\left((\bI_n - \bH) \bT \bQ \bX (\bX^\T \bQ \bX)^{-1} \bR (\bX^\T \bQ \bX)^{-1} \bX^\T \right),
    $$
    where $\bq = \mathrm{diag}(\bQ)$.
\end{proof}

\subsection{Proof of Theorem \ref{thm: adjustment_interpretation}}

\begin{proof}
    Let $\bH$ and $\hat{\bm{\varepsilon}} = (\hat{\varepsilon}_1, \ldots, \hat{\varepsilon}_n)^\T$ be the hat matrix and residual vector of the full model respectively. Write $\hat{\bm{\varepsilon}}_{\rm loocv} = (\hat{\varepsilon}_1^{-1}, \ldots, \hat{\varepsilon}_n^{-n})^\T$. Note that
    $\hat{\bm{\varepsilon}} = (\bI_n - \bH) \by$, $\hat{\bm{\varepsilon}}_{\rm loocv} = \mathrm{diag}((1 - h_{ii})^{-1}) (\bI_n - \bH) \by$, and $\var(\by \vert \bX) = \sigma_\varepsilon^2 \bT$. Then, we have
    $$
        \var(\bW^{\frac{1}{2}} \hat{\bm{\varepsilon}} \vert \bX) = \sigma_\varepsilon^2 \bW^{\frac{1}{2}} (\bI_n - \bH) \bT (\bI_n - \bH)^\T \bW^{\frac{1}{2}},
    $$
    and
    $$
        \var(\bW^{\frac{1}{2}} \hat{\bm{\varepsilon}}_{\rm loocv} \vert \bX) = \sigma_\varepsilon^2 \mathrm{diag}\left(\frac{\sqrt{w_i}}{1 - h_{ii}}\right) (\bI_n - \bH) \bT (\bI_n - \bH)^\T \mathrm{diag}\left(\frac{\sqrt{w_i}}{1 - h_{ii}}\right).
    $$
    Thus,
    $$
    	\begin{aligned}
    		\trace(\bA \bT) 
    		& = \trace[(\bI_n - \bH)^\T \mathrm{diag}(w_i (1-h_{ii})^{-2}) (\bI_n - \bH) \bT] - \trace[(\bI_n - \bH)^\T \mathrm{diag}(w_i) (\bI_n - \bH) \bT]\\
    		& = \frac{1}{\sigma_\varepsilon^2} (\trace[\var(\bW^{\frac{1}{2}} \hat{\bm{\varepsilon}}_{\rm loocv} \vert \bX)] - \trace[\var(\bW^{\frac{1}{2}} \hat{\bm{\varepsilon}} \vert \bX)]) \\
    		& = \frac{1}{\sigma_\varepsilon^2} \sum_{i=1}^n w_i [\var(\hat{\varepsilon}^{-i}_i \vert \bX) - \var(\hat{\varepsilon}_i \vert \bX)].
    	\end{aligned}
    $$
    
    On the other hand, since $\dfR$ is defined through the excess variance of the predictions on the test data against those on the training data, it can be rewritten as
    $$
    	\begin{aligned}
    		\dfR 
    		& = \frac{1}{2 \bar{w} \sigma_\varepsilon^2} \left( n \E[w_\ast \var(\hat{\mu}_\ast \vert \bx_\ast, \bX) \vert \bX] - \sum_{i=1}^n w_i [\var(\hat{\mu}_i \vert \bX) - 2 \cov(y_i, \hat{\mu}_i \vert \bX)]\right)\\
    		& = \frac{1}{2 \bar{w} \sigma_\varepsilon^2} \left( n \E[w_\ast \var( \hat{\varepsilon}_\ast \vert \bx_\ast, \bX) \vert \bX] -  \sum_{i=1}^n w_i \var(\hat{\varepsilon}_i \vert \bX)  - n \sigma_\varepsilon^2 \E(w_\ast \tau_\ast) + \sigma_\varepsilon^2 \sum_{i=1}^n w_i \tau_i\right).
    	\end{aligned}
    $$
    Therefore, we have
    $$
    	\frac{1}{n}\sigma_\varepsilon^2 \xi_\bX = \E[w_\ast \var(\hat{\varepsilon}_\ast \vert \bx_\ast, \bX) \vert \bX] - \frac{1}{n} \sum_{i=1}^n w_i \var(\hat{\varepsilon}^{-i}_i \vert \bX) - \sigma_\varepsilon^2 \left( \E(w_\ast \tau_\ast) - \frac{1}{n}\sum_{i=1}^n w_i \tau_i\right). 
    $$
\end{proof}

\newpage
\section{Data Analysis}

\subsection{Cancer Mortality Data} \label{appendix: rda_cancer}

Table \ref{tab: rda_cancer_coef} shows the variables and their estimated coefficients of the models selected by $\widehat{\wErrR}$ for both the OLS and WLS methods based on a particular partition of the data into the training and test sets. In this example, both methods select the first 10 variables, including the categorical variable \texttt{Region} that consists of 3 categories: Midwest, Northeast and South. We see that the estimated coefficients for these 3 regions have relatively large change when case weights are incorporated into the model. Since the weights are estimated based on the population size, and county population varies a lot across different regions, the change in these coefficients is a direct result of the weight scheme. On the other hand, the contrasts of these coefficients as well as the coefficients of other numerical variables remain relatively stable from the OLS to the WLS model, suggesting that the weight scheme affects more on the variance structure of the model rather than the regression mean. 
\begin{table}[h]
	\centering
	\caption{Estimated coefficients in the optimal models identified by $\widehat{\wErrR}$ for the OLS and WLS methods.}
	\label{tab: rda_cancer_coef}
	\begin{threeparttable}
		\fontsize{8}{9.6}\selectfont
		\begin{tabular}{p{0.2\textwidth}p{0.35\textwidth}>{\centering}p{0.15\textwidth}>{\centering}p{0.08\textwidth}>{\centering\arraybackslash}p{0.08\textwidth}}
			\hline
			Variable & Description & Transformation & OLS & WLS \\
			\hline
			\texttt{RegionMidwest} & Whether the county is in the Midwest & - & $43.222$ \qquad * & $49.778$ \qquad * \\
			
			\texttt{RegionNortheast} & Whether the county is in the Northeast & - & $32.278$ & $38.988$ \qquad * \\
			
			\texttt{RegionSouth} & Whether the county is in the South & - & $51.235$ \qquad * & $57.468$ \qquad *\\
			
			\texttt{avgDeathRateEst2015} & Estimated death rate based on 2015 population estimates and average number of reported cancer mortalities from 2010 to 2016 & - & $0.243$ \quad *** & $0.237$ \quad *** \\
			
			\texttt{incidenceRate} & Mean \textit{per capita} (100,000) cancer diagnoses & - & $0.209$ \quad *** & $0.195$ \quad *** \\ 
			
			\texttt{PctPublicCoverage} & Percent of county residents with government provided health coverage & logit & $-41.652$ *** & $-38.457$  *** \\
			
			\texttt{PctBachDeg25\_Over} & Percent of county residents aged 25 and over with bachelor's degree as the highest education attained & logit & $-9.127$\quad* & $-10.532$ \quad** \\
			
			\texttt{PctMarried} & Percent of county residents who are married & logit & $-9.351$ & $-11.408$ \quad $\cdot$ \\
			
            \texttt{PctPrivateCoverage} & Percent of county residents with private health coverage  & logit &  $-27.979$  *** & $-23.903$  *** \\ 
			
			\texttt{AvgHouseholdSize} & Average number of people in a household & - & $-0.148$ & $-0.833$\\
			
			\texttt{AvgAnnCount} & Average number of reported cases of cancer diagnosed annually & log & $-1.078$ & $-0.851$\\
			
			\texttt{PctUnemployed16\_over} & Unemployment rate for county residents aged 16 and over & logit & $11.463$ *** & $12.204$ ***\\
			\hline
		\end{tabular}
		
		\begin{tablenotes}\footnotesize
			\item \hfill Significance codes: `***' $ < 0.001$, `**' $0.001 \sim 0.01$, `*' $0.01 \sim 0.05$, `.' $0.05 \sim 0.1$, ` ' $>0.1$.
		\end{tablenotes}
	\end{threeparttable}
\end{table}

\subsection{King County House Sales Data} \label{appendix: rda_house_sales}

\begin{figure}[!t]
    \centering
    \includegraphics[scale = 0.52]{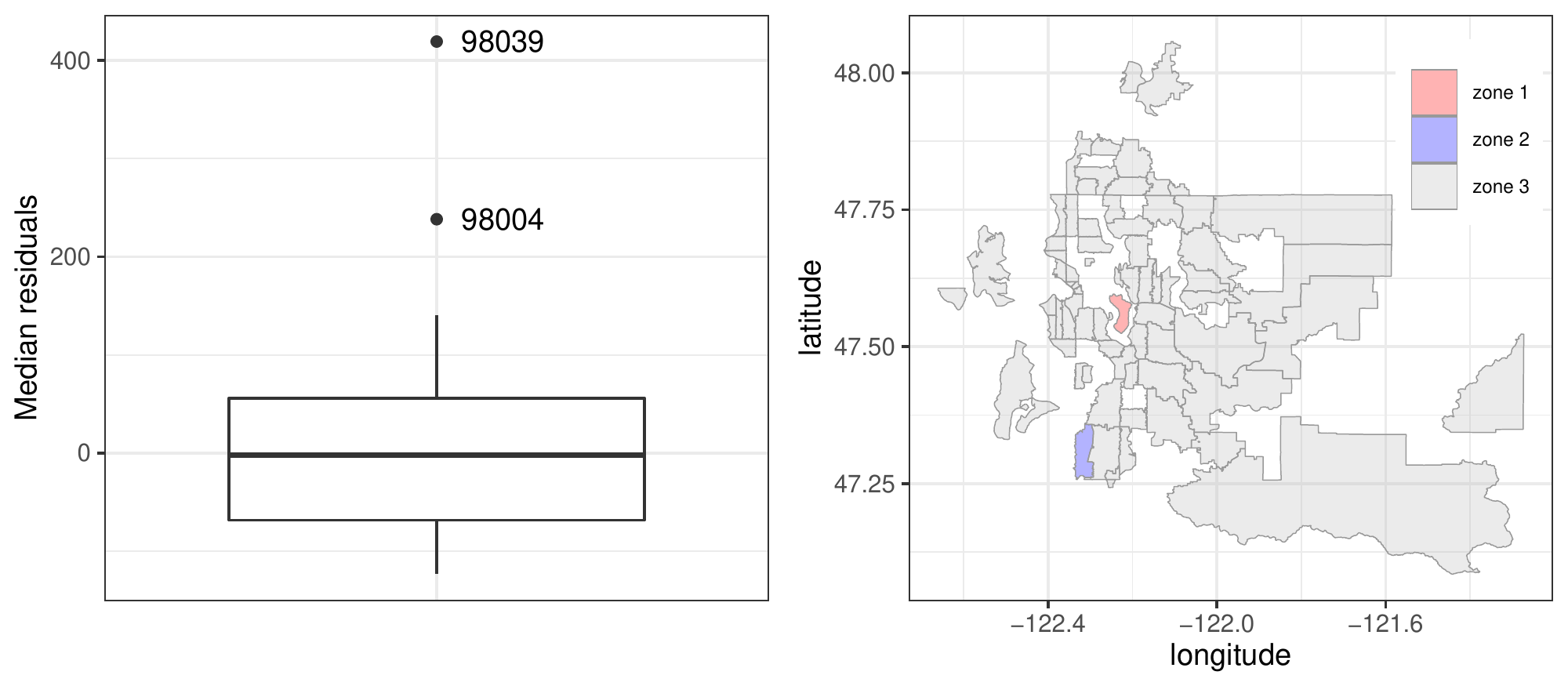}
    \caption{Boxplot of the median residuals of the OLS model for each zip code (left) and the map of the three zones defined through the boxplot (right).}
    \label{fig: rda_house_zone_map}
    \vspace{2em}
    
    \includegraphics[scale = 0.52]{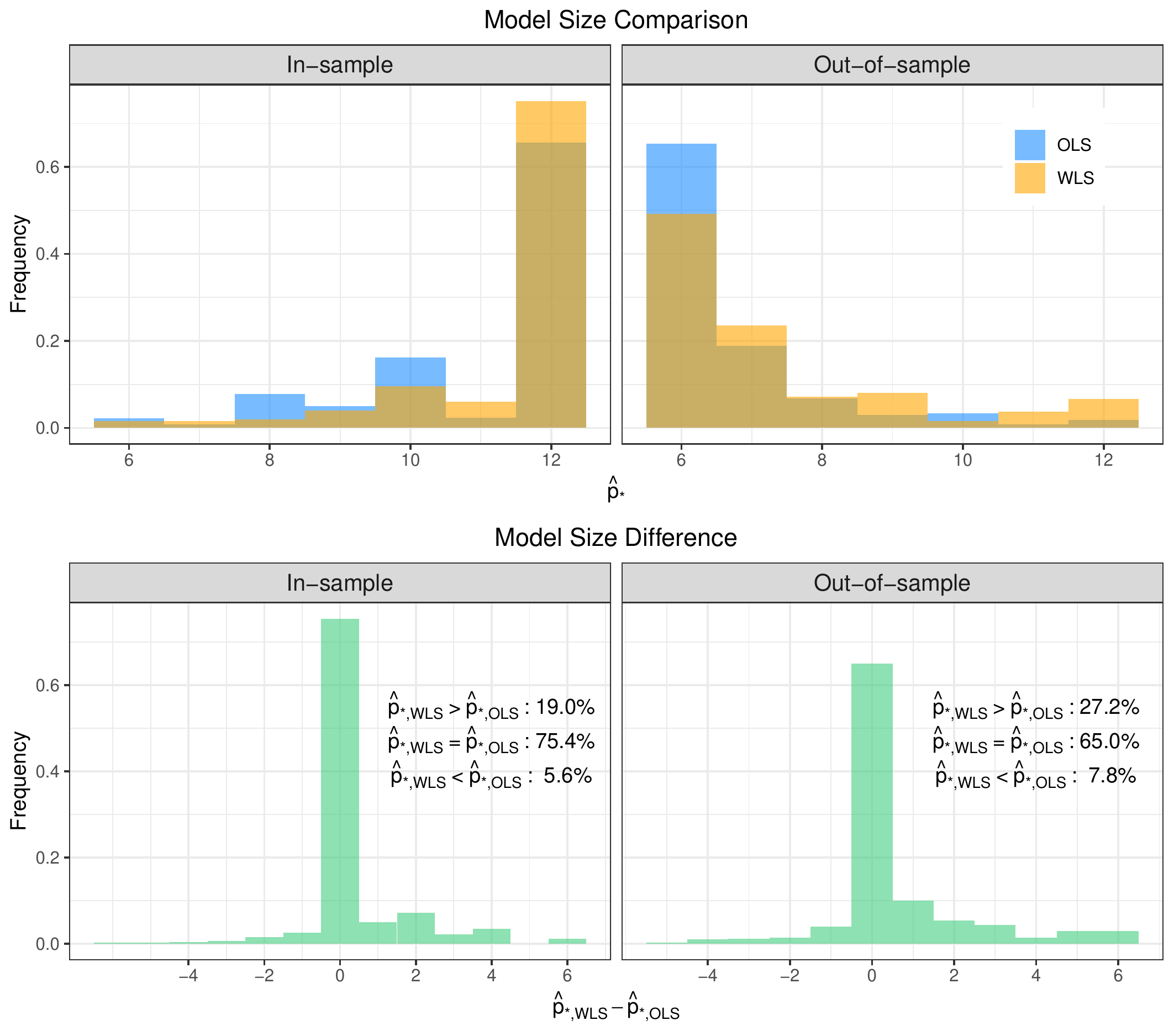}
    \caption{Comparison of the selected model size $\hat{p}_\ast$ between the OLS model and WLS model based on 500 random subsets of size $n = 200$ from the King County house sales data in 2015. The upper panels compare the model size directly whereas the lower panels show the difference in model size.} \label{fig: rda_house_model_selection}
\end{figure}

In this section, we provide more details on the King County housing data analysis described in Section 6.2 of the paper. Table \ref{tab: rda_house_var_list} lists all 12 variables considered in the regression analysis. To simplify the analysis, some of these variables are redefined based on variables in the original dataset, including
\begin{itemize}
    \item \texttt{bedrooms}: We convert the original integer-valued variable to a categorical variable of 3 levels.
    \item \texttt{condition}: We map the original condition indexed from 1 (worst) to 5 (best) to a categorical variable of 3 levels: poor (1 and 2), average (3) and good (4 and 5). 
    \item \texttt{yr\_built}: The variable is initially the year when the house was built. We group its values into 6 time periods.
    \item \texttt{zone}: This variable summarizes \texttt{zipcode} (the zip code of the area in which the house is located) in the original dataset into 3 zones. To determine the grouping, we fit an OLS model with the first 11 variables in Table \ref{tab: rda_house_var_list} using all observations in 2014. The residuals are then grouped by the zip code and the median residuals for each zip code are studied through a boxplot shown in the left panel of Figure \ref{fig: rda_house_zone_map}. Two outlying zip codes with much higher median residuals are identified. We define these two zip codes as zone 1 (98039) and zone 2 (98004) respectively. The remaining zip codes are merged into one common zone called zone 3.  
\end{itemize}
\begin{table}[!b]
	\centering
	\caption{List of independent variables used in King County house sales data analysis}
	\label{tab: rda_house_var_list}
	\fontsize{9}{10.8}\selectfont
	\begin{tabularx}{\textwidth}{lcXc}
		\hline
		Variable & Type & Description & Transformation\\
		\hline
		\texttt{basement} & categorical & Whether the house has basement or not: yes/no  & - \\ 
		
		\texttt{bathrooms} & numeric & Number of bathrooms in the house & - \\ 
		
		\texttt{bedrooms} & categorical & Number of bedrooms in the house: 1-2/3-5/6 or more  & - \\
		
		\texttt{condition} & categorical & Condition of the house: poor/average/good &  \\
		
		\texttt{grade} & numeric & Overall grade of the house regarding the types of materials used and the quality of workmanship & - \\
		
		\texttt{grade2} & numeric & Squared of \texttt{grade} & - \\
		
		\texttt{sqft\_living} & numeric & Interior living space of the houses in square feet & square root \\
		
		\texttt{sqft\_living15} & numeric & Average interior living space for the closest 15 houses in square feet &  square root\\
		
		\texttt{view} & numeric & Grade of the view around the house & -\\
		
		\texttt{waterfront} & categorical & Whether the house has a waterfront or not: yes/no & -\\
		
		\texttt{yr\_built} & categorical & Year when the house was built: 1920 or before/1921-1940/1941-1960/1961-1980/1981-2000/2001 or after & -\\
		
		\texttt{zone} & categorical & Location of the house: zone 1/zone 2/zone 3 & -\\
		\hline
	\end{tabularx}
\end{table}

To estimate the variance function $\tau$, we first fit an OLS model with all variables included using the samples in 2014. The residuals are then bucketed by 100 percentiles of the fitted values, and their variances are calculated and plotted against the median fitted values for each bucket, as shown in Figure 10 in the paper. It turns out that a quadratic function of the fitted value fits the error variances well. This is approximately equivalent to assuming that the variance function $\tau$ is of the form
$$
    \tau(\bx) = a(1 + b \vert \mu \vert )^2, \quad a>0, \,b \neq 0.
$$
With the estimated variance function, the variable sequence is then determined using a similar approach described in Section \ref{subsubsec: rda_cancer} based on all samples in 2014.

\begin{figure}[!t]
    \centering
    \includegraphics[scale = 0.52]{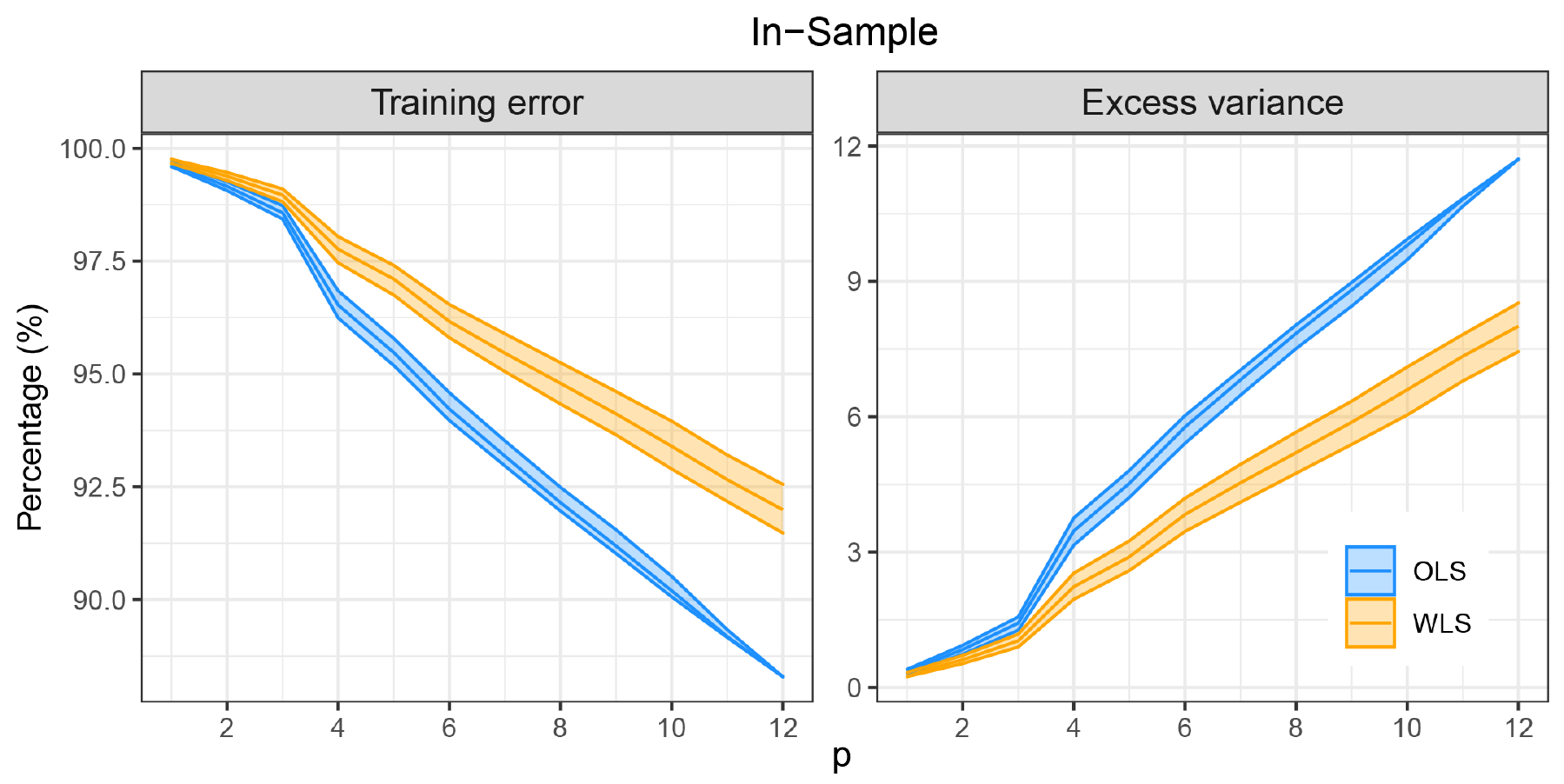}\\
    \vspace{8pt}
    \includegraphics[scale = 0.52]{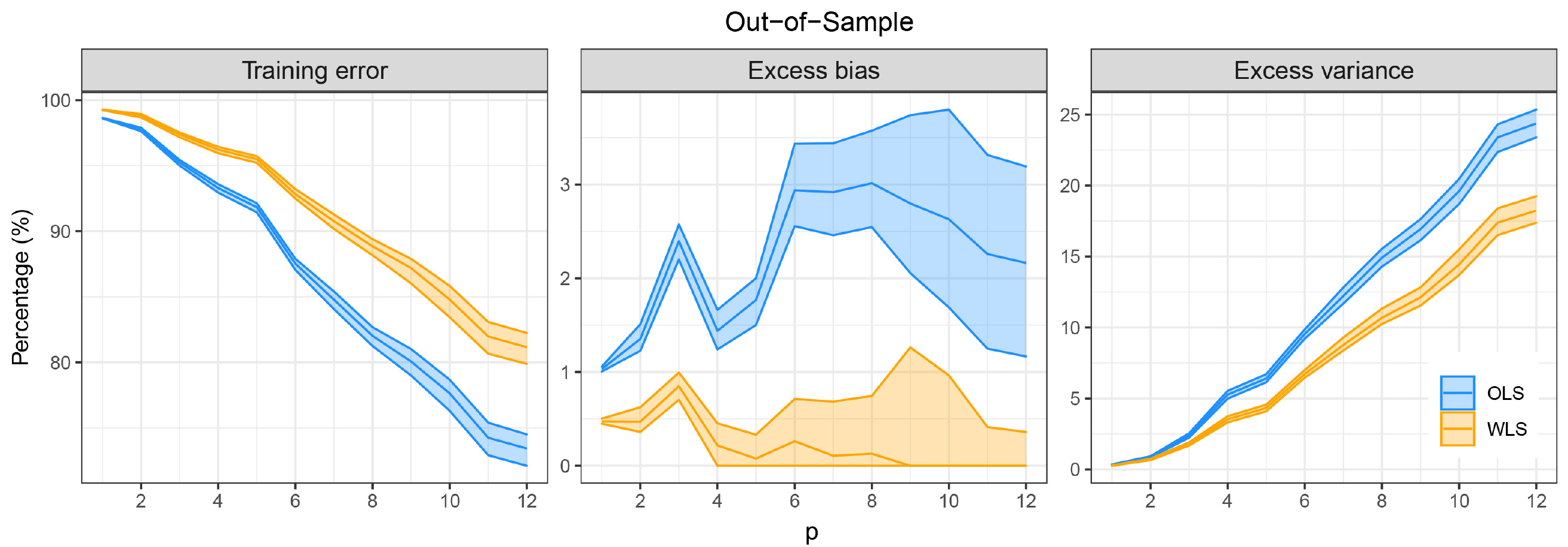}
    \caption{Median percentages of the training error, excess bias and excess variance in the prediction error estimates $\widehat{\wErrF}$ and $\widehat{\wErrR}$ of the selected model as the model size $p$ varies. The shaded bands represent the 2.5th and 97.5th sample percentiles for each $p$.}
    \label{fig: rda_house_component_pct}
\end{figure}

Figure \ref{fig: rda_house_model_selection} shows the size of the selected model for the two procedures based on 500 random training subsets, whereas Figure \ref{fig: rda_house_component_pct} presents the decompositions of $\widehat{\wErrF}$ and $\widehat{\wErrR}$ into the training error, excess bias and excess variance for the same training subsets.

\end{appendices}

\end{document}